\documentclass[10pt]{amsart}
\usepackage{amsmath, amsfonts,amssymb,amscd}

\title[Weak extension theorem for measure-preserving homeomorphisms]
{Weak extension theorem for measure-preserving homeomorphisms of  noncompact manifolds}

\author[T.Yagasaki]{Tatsuhiko Yagasaki}
\address{Division of Mathematics, 
Graduate School of Science and Technology, Kyoto Institute of Technology,  
Matsugasaki, Sakyoku, Kyoto 606-8585, Japan}
\email{yagasaki@kit.ac.jp}

\subjclass[2000]{Primary 57S05; Secondary 58C35}
\keywords{Group of measure preserving homeomorphisms,  
Extension, Principal bundle, Whitney topology, sigma-compact manifolds}

\thanks{This work is supported by Grant-in-Aid for Scientific Research No.19540078.} 

\newtheorem{theorem}{Theorem}[section]
\newtheorem{proposition}{Proposition}[section]
\newtheorem{corollary}{Corollary}[section]
\newtheorem{lemma}{Lemma}[section]

\theoremstyle{definition}
\newtheorem{defi}{Definition}[section] 

\newtheorem{remark}{Remark}[section]
\newtheorem{problem}{Problem}[section]

\newcommand{\lra}{\longrightarrow}

\newcommand{\ds}{\displaystyle} 
\newcommand{\cal}{\mathcal} 
\newcommand{\cbox}{\boxdot} 
\newcommand{\cov}{\mathrm{cov}} 
\newcommand{\lmt}{\longmapsto} 
\def \bs {\boldsymbol}
\def \reg {\text{-reg}}
\def \bireg {\text{-biregular}}
\def \ereg {\text{-e-reg}}
 
\def \e {\varepsilon} 

\newcommand{\llaw[1]}{
\begin{array}[b]{c}
{#1} \\[-1mm] 
\lla
\end{array}}
\newcommand{\lraw[1]}{
\begin{array}[b]{c}
{#1} \\[-1mm] 
\lra
\end{array}}

\newcommand{\U}{\mathcal U}
\newcommand{\V}{\mathcal V}

\begin{document}

\baselineskip 5mm 

\maketitle

\begin{abstract} 
In this paper we deduce weak type extension theorems for the groups of measure-preserving homeomorphisms of  noncompact manifolds. 
As an application, we show that the group of 
measure-preserving homeomorphisms with compact support of a noncompact connected manifold, endowed with the Whitney topology, is locally contractible. 
\end{abstract} 

\section{Introduction}

In this paper we study some topological properties of groups of 
measure preserving homeomorphisms 
and spaces of measure preserving embeddings 
in noncompact manifolds (cf. \cite{Be3, Be4, Fa, Ya1, Ya2}). 
Suppose $M$ is a $\sigma$-compact topological $n$-manifold possibly with boundary and 
$U$ is an open subset of $M$.  
Let ${\mathcal E}^\ast(U, M)$ denote the space of proper embeddings of $U$ into $M$ endowed with the compact-open topology. 
The local deformation lemma for ${\mathcal E}^\ast(U, M)$ \cite{Ce, EK} asserts that 
for any compact subset $C$ of $U$ and any compact neighborhood $K$ of $C$ in $U$ there exists a deformation $\varphi_t$ ($t \in [0,1]$) of 
an open neighborhood ${\mathcal V}$ of the inclusion map $i_U : U \subset M$ in ${\mathcal E}^\ast(U, M)$ such that $\varphi_0(f) = f$, $\varphi_1(f)|_C =  i_C$ and $\varphi_t(f)|_{U - K} =  f|_{U - K}$ ($t \in [0,1]$) for each $f \in {\mathcal V}$. 
For a subset $A$ of $M$ let ${\cal H}_{A}(M)$ denote the group of homeomorphisms $h$ of $M$ with $h|_A = id_A$ endowed with the compact-open topology. 
The local deformation lemma is equivalent to the following weak type extension theorem:  
for any compact neighborhood $L$ of $C$ in $U$ 
there exists a neighborhood ${\cal V}$ of $i_U$ in  ${\cal E}^\ast(U, M)$ and 
a homotopy $s_t : {\cal U} \to {\cal H}_{M - L}(M)$ such that 
$s_0(f) = id_M$ and $s_1(f)|_C = f|_C$ $(f \in {\cal U})$. 

This result motivates the following general formulation: 
Suppose $G$ is a topological group acting on $M$ with the unit element $e$.  Consider the subspace of ${\mathcal E}^\ast(U, M)$ defined by 
${\mathcal E}^G(U, M) = \{ \widehat{g}|_U \mid g \in G\}$, 
where $\widehat{g}$ denotes the homeomorphism on $M$ induced by $g \in G$. 
The weak extension theorem for the group action of $G$ on $M$ asserts that 
there exists a neighborhood ${\cal U}$ of $i_U$ in  ${\cal E}^G(U, M)$ and 
a homotopy $s_t : {\cal U} \to G$ such that 
$s_0(f) = e$ and $\widehat{s_1(f)}|_C = f|_C$ $(f \in {\cal U})$. 

Suppose $\mu$ is a good Radon measure on $M$ with $\mu(\partial M) = 0$.  
Let ${\cal H}(M; \mu)$ and ${\cal H}(M; \mu\reg)$ denote 
the subgroups of ${\cal H}(M)$ consisting of 
$\mu$-preserving homeomorphisms and $\mu$-biregular homeomorphisms of $M$ and let 
${\cal E}^\ast(U, M; \mu\reg)$ denote the subspace of 
${\cal E}^\ast(U, M)$ consisting of 
$\mu$-biregular proper embeddings of $U$ into $M$. 
In \cite{Fa} A. Fathi obtained a local deformation lemma for the space ${\cal E}^\ast(U, M; \mu\reg)$ (\cite[Theorem 4.1]{Fa}). 
This is reformulated as the weak extension theorem for the group ${\cal H}(M; \mu\reg)$ (\cite[Corollary 4.2]{Fa}). 
In the case $M$ is compact and connected, 
he also obtained a selection theorem for $\mu$-biregular measures on $M$ (\cite[Theorem 3.3]{Fa}) and 
used these results to deduce the weak extension theorem for the group ${\cal H}(M; \mu)$ (\cite[Theorem 4.12]{Fa}). 

In this paper we are concerned with the case where $M$ is non-compact. 
In \cite{Be3} R. Belanga has already extended the selection theorem for $\mu$-biregular measures 
to the non-compact case (\cite[Theorem 4.1]{Be3}).  
We combine these results to obtain the weak extension theorem for 
the group ${\cal H}(M; \mu)$ (cf.\,Corollary~\ref{cor_ext_mp}). 

\begin{theorem}\label{main_thm_mp} 
Suppose $M$ is an $n$-manifold, 
$\mu$ is a good Radon measure on $M$ with $\mu(\partial M) = 0$, 
$C$ is a compact subset of $M$, $U$ is an open neighborhood of $C$ in $M$. 
Then there exists a neighborhood ${\cal U}$ of $i_U$ in ${\cal E}^{{\cal H}(M; \mu)}(U, M)$ and a homotopy $s :  {\cal U} \times [0,1] \to {\cal H}(M; \mu)$ such that
\begin{itemize}
\item[(1)] for each $f \in {\cal U}$ \\[0.5mm]
{\rm (i)} $s_0(f) = id_M$, \ {\rm (ii)} 
$s_1(f)|_C = f|_C$, \ {\rm (iii)} 
if $f = id$ on $U \cap \partial M$, then $s_t(f) = id$ on $\partial M$ $(t \in [0,1])$, 
\vskip 0.5mm 
\item[{\rm (2)}] $s_t(i_U) = id_M$ $(t \in [0,1])$. 
\end{itemize} 
\end{theorem} 

In comparison with topological or $\mu$-biregular homeomorphisms, 
``$\mu$-preserving homeomorphism'' is a global property and 
we can not obtain a compactly supported weak extension theorem for the group ${\cal H}(M; \mu)$. This obstruction vanishes on the kernel of the end charge homomorphism $c^\mu$. 

In \cite{AP2} S.~R.~Alpern and V.\,S.\,Prasad introduced the end charge homomorphism $c^\mu$, which is a continuous homomorphism 
defined on the subgroup ${\cal H}_{E_M}(M; \mu)$ 
of $\mu$-preserving homeomorphisms of $M$ which fix the ends of $M$. 
The kernel of $c^\mu$, ${\rm ker}\,c^\mu$, includes the subgroup 
${\cal H}_c(M; \mu)$ of $\mu$-preserving homeomorphisms of $M$ with compact support. 
If $h \in {\cal H}_{E_M}(M, E; \mu)$ and $c^\mu(h) = 0$, then 
one can split moves of $\mu$-volume by $h$. 
Hence, we can obtain the compactly supported weak extension theorem for the subgroup ${\rm ker}\,c^\mu$ (cf.\,Theorem~\ref{thm_ext_ker}). 

\begin{theorem}\label{main_thm_ext_ker} 
Suppose $M$ is a connected $n$-manifold, 
$\mu$ is a good Radon measure on $M$ with $\mu(\partial M) = 0$, 
$C$ is a compact subset of $M$ and 
$U$ and $V$ are open neighborhoods of $C$ in $M$ such that  
$V \cap O$ is connected for each connected component $O$ of $M - C$. 
Then there exists a neighborhood ${\cal U}$ of $i_U$ in ${\cal E}^{{\rm ker}\,c^\mu}(U, M)$ and a homotopy $s :  {\cal U} \times [0,1] \to {\cal H}_{M  - V, c}(M; \mu)$ such that
\begin{itemize}
\item[(1)] for each $f \in {\cal U}$ \\[0.5mm]
{\rm (i)} $s_0(f) = id_M$, \ {\rm (ii)} 
$s_1(f)|_C = f|_C$, \ {\rm (iii)} 
if $f = id$ on $U \cap \partial M$, then $s_t(f) = id$ on $\partial M$ $(t \in [0,1])$, 
\vskip 0.5mm 
\item[{\rm (2)}] $s_t(i_U) = id_M$ $(t \in [0,1])$. 
\end{itemize} 
\end{theorem} 

We also discuss a non-ambient deformation lemma for $\mu$-preserving embeddings (Theorem~\ref{thm_non-amb_def}). 

In the last section we study the group ${\cal H}_c(M; \mu)_w$ endowed with the Whitney topology (cf.\,\cite{BMSY}). It is known that the group
${\cal H}(N)$ and the subgroup 
${\cal H}(N; \nu)$ are locally contractible 
for any compact $n$-manifold $N$ and any good Radon measure $\nu$ on $N$ with $\nu(\partial N) = 0$ 
(\cite[Corollary 1.1]{EK}, \cite[Theorem 4.4]{Fa}). 
In \cite{BMSY} it is shown that the group ${\cal H}_c(M)_w$ consisting of homeomorphisms of $M$ with compact support, endowed with the Whitney topology, is locally contractible. 
In this article, as an application of the weak extension theorem for ${\cal H}_c(M; \mu)$, 
we show that the group ${\cal H}_c(M; \mu)_w$ is also locally contractible 
for any connected $n$-manifold $M$ (Theorem~\ref{thm_LC}). 

This paper is organized as follows. 
Section 2 is devoted to the general formulations and basic properties 
of local weak extension property and local weak section property for group actions. 
Section 3 contains fundamental facts related to Radon measures on manifolds 
(selection theorems for measures, end charge homomorphism, etc.). 
In Section 4 we recall  
the local deformation lemma for biregular embeddings and 
discuss some direct consequences of this lemma. 
In Section 5 we obtain   
the weak extension theorems for the groups ${\cal H}(M; \mu)$, ${\rm ker}\,c^\mu$ and 
${\cal H}_c(M; \mu)$ and a non-ambient deformation lemma for $\mu$-preserving embeddings. 
In Section 6 we recall basic facts on the Whitney topology and 
show that the group ${\cal H}_c(M; \mu) _w$ is 
locally contractible for any connected $n$-manifold $M$. 

\section{Fundamental facts on group actions}  

\subsection{Conventions} \mbox{} 

For a topological space $X$ and a subset $A$ of $X$, 
the symbols ${\rm Int}_X A$,  $cl_X A$ and ${\rm Fr}_X A$ denote the topological 
interior, closure and frontier of $A$ in $X$. 
Let ${\cal C}(X)$ denote the collection of all connected components of $X$. 

Suppose $Y$ is a locally connected, locally compact Hausdorff space. 
Let ${\cal H}(Y)$ denote the group of homeomorphisms of $Y$ 
endowed with the compact-open topology. 
For a subset $A$ of $Y$, let  
${\cal H}_A(Y) = \big\{h \in {\cal H}(Y) \mid h|_A = id_A \big\}$ (with the subspace topology). 
The group ${\cal H}(Y)$ and the subgroup ${\cal H}_A(Y)$ are topological groups. 
In general, for any topological group $G$, 
the symbols $G_0$ and $G_1$ denote the connected component and the path-component of the unit element $e$ in $G$. 

For subspaces $A \subset X$ of $Y$ let ${\cal E}(X, Y)$ denote the space of embeddings 
$f : X \hookrightarrow Y$ endowed with the compact-open topology, and 
let ${\cal E}_A(X, Y) = \big\{ f \in {\cal E}(X, Y) \mid f|_A = id_A \big\}$ (with the subspace topology). 
By $i_X : X \subset Y$ we denote the inclusion map of $X$ into $Y$.  

In this article, an $n$-manifold means a paracompact $\sigma$-compact (separable metrizable) topological $n$-manifold {\em possibly with boundary}. 
Suppose $M$ is an $n$-manifold. 
The symbols $\partial M$ and ${\rm Int}\,M$ denote the boundary and interior of $M$ as a manfiold. 
For a subspace $X$ of $M$, 
an embedding $f : X \to M$ is said to be {\em proper} if $f^{-1}(\partial M) = X \cap \partial M$. 
Let ${\cal E}^\ast(X, M)$ denote the subspace of ${\cal E}(X, M)$ consisting of 
proper embeddings $f : X \to M$. 
For a subset $A$ of $X$ let ${\cal E}_A^\ast(X, M) = {\cal E}^\ast(X, M) \cap {\cal E}_A(X, M)$. 

By an {\em $n$-submanifold} of $M$ we mean 
a closed subset $N$ of $M$ such that $N$ is an $n$-manifold and 
${\rm Fr}_MN$ is locally flat in $M$ and transverse to $\partial M$ so that 
(i) $M - {\rm Int}_M N$ is an $n$-manifold and (ii) 
${\rm Fr}_M N$ and $N \cap \partial M$ are $(n-1)$-manifolds with the common boundary $({\rm Fr}_M N) \cap (N \cap \partial M)$.  
For simplicity, let $\partial_+ N = {\rm Fr}_M N$, $\partial_- N = N \cap \partial M$ and $N^c = M - {\rm Int}_M N$. 
More generally, for a subset $U$ of $M$ let $\partial_- U = U \cap \partial M$. 

Suppose $M$ is an $n$-manifold. 

\begin{lemma}\label{lem_submfd} {\rm (\cite[Theorem 0]{AP1}, cf.\,\cite{KS})}  
Suppose $C$ is a compact subset of $M$ and $U$ is a neighborhood of $C$ in $M$. Then there exists a compact $n$-submanifold $N$ of $M$ such that 
$C \subset {\rm Int}_M N$ and $N \subset U$. 
\end{lemma} 

\begin{lemma}\label{lem_connect} 
\begin{itemize}
\item[] \hspace{-13mm} {\rm (1)} If $M$ is connected and 
$L$ is an $n$-submanifold of $M$ such that $\partial_+ L$ is compact, then there exists a connected $n$-submanifold $N$ of $M$ such that 
$L \subset {\rm Int}_M N$ and $N \cap L^c$ is compact. 
\item[(2)] Suppose $C$ is a compact subset of $M$. 
\begin{itemize}
\item[(i)\,] For any neighborhood $U$ of $C$ in $M$ there exists a compact $n$-submanifold $N$ of $M$ such that 
$C \subset {\rm Int}_M N$, $N \subset U$ and $O - N$ is connected for each $O \in {\cal C}(M - C)$. 
\item[(ii)] If $U$ is an open neighborhood of $C$ in $M$ such that 
$U \cap O$ is connected for each $O \in {\cal C}(M - C)$, then 
there exists a compact $n$-submanifold $N$ of $M$ such that 
$C \subset {\rm Int}_M N$, $N \subset U$ and 
$N \cap O$ is connected for each $O \in {\cal C}(M - C)$.
\end{itemize}
\end{itemize}
\end{lemma} 

\begin{proof} 
(1) Since $M$ is connected and $\partial_+ L$ is compact, 
${\cal C}(L)$ is a finite collection. 
Since $M$ is connected, there exists 
a finite collection of disjoint arcs $\{ \alpha_i \}_i$ in $L^c$ such that 
$L \cup (\bigcup_i \alpha_i)$ is connected. 
We apply Lemma~\ref{lem_submfd} to 
$C = \partial_+ L \cup (\bigcup_i \alpha_i)$ in the $n$-manifold $L^c$ 
in order to find a compact $n$-submanifold $N_0$ of $L^c$ such that  
$C \subset {\rm Int}_{L^c} N_0$ and each $K \in {\cal C}(N_0)$ meets $C$. 
Then $N = L \cup N_0$ satisfies the required conditions. 
\vskip 1mm 
(2)\,(i) We may assume that $M$ is connected 
(apply the connected case to each component of $M$). 
By Lemma~\ref{lem_submfd} there exists a compact $n$-submanifold $N_1$ of $M$ such that 
$C \subset {\rm Int}_M N_1$ and $N_1\subset U $. 
Let ${\cal C} = \big\{ O \in {\cal C}(M - C) \mid O \not \subset N_1 \}$. 
Since ${\cal C}(N_1^c)$ is a finite collection, so is ${\cal C}$. 

For each $O \in {\cal C}$, it is seen that 
$O$ is a connected $n$-manifold, 
$N_1^c \cap O$ is an $n$-submanifold of $O$, 
$(N_1^c \cap O)^c = N_1 \cap O$ in $O$ and 
${\rm Fr}_O (N_1^c \cap O) = ({\rm Fr}_M N_1) \cap O$ is compact 
(it is a union of components of ${\rm Fr}_M N_1$). 
Thus, by (1) we can find a connected $n$-submanifold $L_O$ of $O$ such that 
$N_1^c \cap O \subset {\rm Int}_O L_O$ and $L_O \cap (N_1 \cap O)$ is compact. 
Note that $L_O$ is closed in $M$ so that it is also a connected $n$-submanifold of $M$. 
Let $L = \bigcup_{O \in {\cal C}} L_O$. Then, $N = L^c$ satisfies the required conditions. 
In fact, $C \subset M - L = {\rm Int}_M N$, $N \subset N_1$, 
${\cal C} =  \big\{ O \in {\cal C}(M - C) \mid O \not \subset N \}$ 
and $O - N = {\rm Int}_M L_O$ for each $O \in {\cal C}$. 

(ii) Since ${\cal C}(U - C) = \{ O \cap U \mid O \in {\cal C}(M - C) \}$,  
by replacing $M$ by $U$, we may assume that $U = M$. 
Again we may assume that $M$ is connected. 
By Lemma~\ref{lem_submfd} there exists a compact $n$-submanifold $N_1$ of $M$ such that 
$C \subset {\rm Int}_M N_1$. 
Consider the finite collection ${\cal C} = \big\{ O \in {\cal C}(M - C) \mid O \not \subset N_1 \}$. 
For each $O \in {\cal C}$, it is seen that 
$O$ is a connected $n$-manifold, 
$N_1 \cap O$ is an $n$-submanifold of $ O$, 
$(N_1 \cap O)^c = N_1^c \cap O$ in $O$ and 
${\rm Fr}_{O} (N_1 \cap O) = ({\rm Fr}_M N_1) \cap O$ is compact. 
Thus, by (1) we can find a connected $n$-submanifold $K_O$ of $O$ such that 
$N_1 \cap O \subset {\rm Int}_{O} K_O$ and $K_O \cap (N_1^c \cap O)$ is compact. 
Then, $N = N_1 \cup (\bigcup_{O \in {\cal C}} K_O)$ satisfies the required conditions. 
In fact, 
$\big\{ O \in {\cal C}(M - C) \mid O \not \subset N \} \subset {\cal C}$ 
and $N \cap O = K_O$ for each $O \in {\cal C}$. 
\end{proof} 

\subsection{Pull-backs} \mbox{} 

For maps $B_1 \lraw[p] B \llaw[\pi]  E$, 
we obtain 
the {\em pull-back} diagram in the category of topological spaces and continuous maps : 
\vspace{-2mm} 
$$\begin{array}[c]{cccl}
& p' & & \\[-0.5mm] 
p^\ast E & \lra & E & \\[1mm] 
\left\downarrow \makebox(0,9){} \right. & \hspace{-25mm} \pi' & \left\downarrow \makebox(0,9){} \right.  & \hspace{-3mm} \pi \\[1mm] 
B_1 & \lra & B & \\[-1mm] 
& p & 
\end{array}$$ 
\noindent 
Explicitly, 
the space $p^\ast E$ and 
the maps $B_1 \llaw[\pi'] p^\ast E \lraw[p']  E$ are defined by 
\[ \mbox{
$p^\ast E = \{ (b_1, e) \in B_1 \times E \mid p(b_1) = \pi(e) \}$ \hspace{2mm}  and \hspace{2mm}  
$\pi'(b_1, e) = b_1$, \ \ $p'(b_1, e) = e.$
}\] 

Suppose a topological group $G$ acts on spaces $B$ and $B_1$ transitively. 
Let $p : B_1 \to B$ be a $G$-equivariant map. 
Fix a point $b_1 \in B_1$ and let $b = p(b_1) \in B$ and 
let $G_b$ be the stabilizer of $b$ under the $G$-action on $B$. 
Consider the orbit map $\pi : G \to B$, $\pi(g) = gb$. 
\vspace{1mm} 
Then the maps $B_1 \lraw[p] B \llaw[\pi]  G$ induce  
the pull-back diagram : \hspace{20mm} 
$\begin{array}[t]{cccl}
& p' & & \\[-0.5mm] 
p^\ast G & \lra & G & \\[1mm] 
\left\downarrow \makebox(0,9){} \right. & \hspace{-25mm} \pi' & \left\downarrow \makebox(0,9){} \right.  & \hspace{-3mm} \pi \\[1mm] 
B_1 & \lra & B & \\[-1mm] 
& p & 
\end{array}$ \\[1mm] 
The group $G_b$ acts freely on $p^\ast G$ on the right by $(x,g) \cdot h = (x, gh)$ ($(x,g) \in p^\ast G$, $h \in G_b$). 
The induced map $p' : p^\ast G \to G$ admits a right inverse 
$r : G \to p^\ast G$, $r(g) = (gb_1, g)$ (i.e., $p'r  = id_G$). 

\begin{defi} We say that the $G$-equivariant map $p : B_1 \to B$ has 
the {\em local section property} for $G$ (LSP$_G$) at $b_1$ 
if there exists a neighborhood $U_1$ of $b_1$ in $B_1$ and 
a map $s_1 : U_1 \to G$ such that $\pi s_1 = p|_{U_1}$. 
\end{defi}

\begin{lemma}\label{lem_p-bdl} 
{\rm (1)} The map $p$ has LSP$_G$ at $b_1$ iff 
the induced map $\pi' : p^\ast G \to B_1$ is a principal $G_b$-bundle. 

{\rm (2)} If the fiber $p^{-1}(b)$ is contractible, then 
the map $p' : p^\ast G \to G$ is a homotopy equivalence. 
\end{lemma}

\begin{proof}
(1) Suppose the map $p$ has LSP$_G$ at $b_1$. 
Take any point $b_2 \in B_1$. 
Since $G$ acts on $B_1$ transitively, there exists a $g \in G$ with $b_2 = gb_1$. 
Then $U_2 = gU_1$ is a neighborhood of $b_2$ in $B_1$ and the map 
$s_2 : U_2 \to G$, $s_2(x) = gs_1(g^{-1}x)$ satisfies 
the condition $\pi s_2 = p|_{U_2}$ 
(i.e., $\pi s_2(x) =  gs_1(g^{-1}x)b = g(p(g^{-1}x)) = p(x)$). 
The map $\pi' : p^\ast G \to B_1$ admits a local trivialization  
$$\phi : U_2 \times G_b \cong (\pi')^{-1}(U_2) = \bigcup_{x \in U_2} \big(\{ x \} \times \pi^{-1}(p(x))\big) \hspace{4mm} \text{over $U_2$ defined by} \hspace{4mm} \phi(x, h) = (x, s_2(x)h).$$ 

The converse is obvious. 

(2) It remains to show that $rp' \simeq id_{p^\ast G}$. 
There exists a contraction $\phi_t : p^{-1}(b) \to p^{-1}(b)$ $(t \in [0,1])$  
such that $\phi_1(p^{-1}(b)) = \{ b_1 \}$. 
If $(x, g) \in p^\ast G$, then $x \in p^{-1}(gb) = g \,p^{-1}(b)$. 
Thus, we can define a homotopy 
\[ \mbox{$\Phi_t : p^\ast G \to p^\ast G$ \hspace{4mm} from \ $id_{p^\ast G }$ \ to \ $rp'$ \ by \hspace{4mm} 
$\Phi_t(x, g) = \big(g \phi_t(g^{-1}x), g)$.} \]
\vskip -7.5mm 
\end{proof} 

\subsection{Group actions and spaces of embeddings} \mbox{} 

Suppose a topological group $G$ acts continuously on a locally compact Hausdorff space $Y$. 
Each $g \in G$ induces $\hat{g} \in {\cal H}(Y)$ defined by $\hat{g}(y) = g y$ $(y \in Y$).  
Let $H$ be any subset of $G$.   
For subsets $A, B$ of $Y$ we have the following 
subsets of $H$: 
$$\begin{array}[c]{l}
H_A = \big\{ h \in H \mid \hat{h}|_A = id_A \big\}, \ \ \ H(B) = H_{Y \setminus B}, \ \ \ 
H_A(B) = H_A \cap H(B), \\[1.5mm] 
H_c = \big\{ h \in H \mid {\rm supp}\,\hat{h} \text{ is compact} \big\}. 
\end{array}$$
If $H$ is a subgroup of $G$, then these are subgroups of $H$. 

For subsets $X \subset C \subset U$ of $Y$, 
the group $G_X(U)$ acts continuously on the space ${\cal E}_X(C, U)$ by the left composition 
$g \cdot f = \hat{g} f$ ($g \in G_X(U)$, $f \in {\cal E}_X(C, U)$) and 
we have the following subspace of ${\cal E}_X(C, U)$: 
$${\cal E}^H_X(C, U) = H_X(U)i_C  = \{ \hat{g}|_C \mid g \in H_X(U) \} 
\quad \text{(with the compact-open topology)}.$$
Since ${\cal E}^H_X(C, U) = {\cal E}^{H_X}(C, U)$, by replacing $H$ by $H_X$ if necessary, we omit $X$ in the subsequent statements. 

Consider the pull-back diagram : 
$$\hspace{20mm} 
\begin{array}[c]{cccl}
& p' & & \\
p^\ast G & \lra & G \\[0.5mm] 
\pi' \left\downarrow \makebox(0,10){} \right. \ \ & & \ \ \left\downarrow \makebox(0,10){}\right. \pi & \hspace{3mm} \text{,where} \ \ \ \pi(g) = \hat{g}|_C \ \ \ \text{and} \ \ \ p(f) = f|_C. \\[2mm] 
{\cal E}^G(U, Y) & \lra & {\cal E}^G(C, Y)& \\ 
& p & & 
\end{array}$$ 
\vskip 1mm 
\noindent The group $G$ acts on the spaces ${\cal E}^G(U, Y)$ and ${\cal E}^G(C, Y)$ transitively. The restriction map $p$ is $G$-equivariant and has the fiber 
$p^{-1}(i_C) = {\cal E}_C^G(U, Y)$. 

\begin{defi} We say that the pair $(U, C)$ has the local section property for $G$ (LSP$_G$) if 
the $G$-equivariant map $p : {\cal E}^G(U, Y) \to {\cal E}^G(C, Y)$ has LSP$_G$
at $i_U$. 
\end{defi}

\begin{lemma}\label{lem_p-bdl_2} The pair $(U, C)$ has LSP$_G$ iff 
the map $\pi' : p^\ast G \to {\cal E}^G(U, Y)$ is a principal $G_C$-bundle.
\end{lemma} 

This lemma follows directly from Lemma~\ref{lem_p-bdl}\,(1).  

\begin{lemma}\label{lem_he} \mbox{} 
Suppose there exists a path $h : [0,1] \to G$ such that $h_0 = e$,  $\widehat{h_1}(U) \subset C$ and $\widehat{h_t}(U) \subset U$, $\widehat{h_t}(C) \subset C$ $(t \in [0,1])$. 
Then the following hold. 
\begin{itemize}
\item[(1)] The map $p : {\cal E}^G(U, Y) \to {\cal E}^G(C, Y)$ is a homotopy equivalence. 

\item[(2)] There exists a strong deformation retraction $\chi_t$ $(t \in [0,1])$ of 
${\cal E}_C^G(U, Y)$ onto the singleton $\{ i_{U} \}$. 

\item[(3)]The map $p' :  p^\ast G \to G$ is a homotopy equivalence. 
\end{itemize}
\end{lemma} 

\begin{proof} 
(1) We can define a map $p_1 : {\cal E}^G(C, Y) \to {\cal E}^G(U, Y)$ by 
 $p_1(f) = f \,\widehat{h_1}|_{U}$.  
It follows that \\[1.5mm] 
\hspace{15mm} 
\begin{tabular}[c]{cl}
(i) & $p_1p(f) = f \,\widehat{h_1}|_{U}$ and 
a homotopy $\phi_t : id \simeq p_1p$ is defined by $\phi_t(f) = f\,\widehat{h_t}|_{U}$, and \\[1.5mm] 
(ii) & $pp_1(f) = f\,\widehat{h_1}|_{C}$ and 
a homotopy $\psi_t : id \simeq pp_1$ is defined by $\psi_t(f) = f\,\widehat{h_t}|_{C}$. 
\end{tabular}
\vskip 1mm
\begin{itemize}
\item[(2)] The contraction $\chi_t$ of ${\cal E}_C^G(U, Y)$ is defined by $\chi_t(f) = \widehat{h_t}^{-1} f \widehat{h_t}|_{U}$. 

\item[(3)] The assertion follows from (2) and Lemma~\ref{lem_p-bdl}\,(2). 
\end{itemize}
\vskip -6mm 
\end{proof}

Lemmas~\ref{lem_p-bdl_2} and \ref{lem_he} yield the following consequence.

\begin{proposition}\label{prop_fib} 
If a subset $C$ of $Y$ satisfies the condition $(\ast)$ below, then the  map 
$$G_C \ \subset \ G \lraw[\pi] {\cal E}^G(C, Y) \hspace{4mm} \text{defined by} \ \  
\pi(h) = \hat{h}|_{C}$$ 
is a locally trivial bundle up to homotopy equivalences and hence has the exact sequence for  homotopy groups. 

\begin{itemize}
\item[$(\ast)$] There exists a subset $U$ of $Y$ such that 
{\rm (i)} $C \subset U$, \ \ {\rm (ii)} the pair $(U, C)$ has LSP$_G$, \ and 
\begin{itemize}
\vskip 1mm
\item[(iii)] there exists a path $h_t \in G$ $(t \in [0,1])$ such that \\[1mm] 
\hspace{10mm} $h_0 = e$, \ \ $\widehat{h_1}(U) = C$, \ \ $\widehat{h_t}(U) \subset U$, \ \ $\widehat{h_t}(C) \subset C$ \ \ $(t \in [0,1])$.
\end{itemize}
\end{itemize}
\end{proposition} 

\subsection{Weak extension property} \mbox{} 

Suppose a topological group $G$ acts on an $n$-manifold $M$. 
Consider a pair $(H, F)$ of subsets of $G$ and a triple $(V, U, C)$ of subsets of $M$ such that $C \subset U \cap V$ (we do not assume that $F \subset H$ and $U \subset V$).  

\begin{defi} 
We say that the triple $(V, U, C)$ has the {\em weak extension property} for $(H, F)$ (abbreviated as WEP$_{H, F}$ or WEP$(H, F)$) if 
there exists a neighborhood ${\cal U}$ of $i_U$ in ${\cal E}^H(U, M)$ and 
a homotopy \break $s :  {\cal U} \times [0,1] \to F(V)$ such that
\begin{itemize}
\item[(1)] for each $f \in {\cal U}$ \\
\hspace{0.5mm} 
{\rm (i)} $s_0(f) = e$, \hspace{0.5mm} {\rm (ii)} 
$\widehat{s_1(f)}|_C = f|_C$, \hspace{0.5mm} {\rm (iii)} 
if $f = id$ on $\partial_- U$, then $\widehat{s_t(f)} = id$ on $\partial M$ $(t \in [0,1])$, 
\item[{\rm (2)}] $s_t(i_U) = e$ $(t \in [0,1])$. 
\end{itemize} 
The map $s_t :  {\cal U} \to F(V)$ ($t \in [0,1]$) is called the 
{\em local weak extension map} (LWE map). 
When $H = F$, we simply say that $(V, U, C)$ has WEP$_H$. 
When $V = U$, we say that the pair $(U, C)$ has WEP$_{H, F}$. 
Note that WEP$_G$ for $(U, C)$ implies LSP$_G$ for $(U, C)$. 
\end{defi} 

One of our interest is the following problem. 

\begin{problem}\label{prob_large}
Given a class of triples $(V, U, C)$ in $Y$ and a subset $F$ of $G$, determine the largest subset $H$ of $G$ for which each triple $(V, U, C)$ in this class has $WEP(H, F)$. 
\end{problem}

The next lemma easily follows from the definition.  

\begin{lemma}\label{l_WEP}
Suppose $(V, U, C)$ and $(V', U', C')$ are two triples of subsets in $M$ such that 
$C \subset U \cap V$ and $C' \subset U' \cap V'$ 
and $(H, F)$ and $(H', F')$ are two pairs of subsets in $G$. 
If {\rm (i)} $(V, U, C)$ has WEP$(H, F)$, 
{\rm (ii)} $V \subset V'$, $U \subset U'$, $C \supset C'$ and {\rm (iii)}
$H \supset H'$, $F \subset F'$, then 
$(V', U', C')$ has WEP$(H', F')$. 
\end{lemma} 

\begin{lemma}\label{l_WEP_disjoint} Suppose $F$ is a subgroup of $G$. 
If two triples $(V_1, U_1, C_1)$ and $(V_2, U_2, C_2)$ have WEP$(H, F)$ and 
$V_1 \cap V_2 = \emptyset$, then the triple 
$\big(V_1 \cup V_2, \ U_1 \cup U_2, \ C_1 \cup C_2 \big)$ also has WEP$(H, F)$. 
\end{lemma} 

\begin{proof} 
For $i = 1, 2$\, let 
$\begin{array}[b]{ccc}
& s^i_t & \\[-0.5mm] 
{\cal E}^H(U_i, M) \ \supset \ {\cal U}_i & \to & F(V_i)
\end{array}$ 
be the associated LWE map for $(V_i, U_i, C_i)$. 
Take a neighborhood ${\cal U}$ of $i_{U_1 \cup U_2}$ in ${\cal E}^H(U_1 \cup U_2, M)$ such that $f|_{U_i} \in {\cal U}_i$ ($i = 1, 2$) for each $f \in {\cal U}$. 
Then the required LWE map 
$s_t : {\cal U} \lra F(V_1 \cup V_2)$
for $\big(V_1 \cup V_2, \ U_1 \cup U_2, \ C_1 \cup C_2 \big)$ is defined by 
$$s_t(f) = s^1_t(f|_{U_1}) s^2_t(f|_{U_2}) \hspace{3mm} (\text{the multiplication in $G$}).$$ 
Note that $\widehat{s_t(f)} = \widehat{s^i_t(f|_{U_i})}$ on $V_i$ and 
$\widehat{s_t(f)} = id$ on $M - (V_1 \cup V_2)$. 
\end{proof}


\section{Spaces of Radon measures and 
groups of measure-preserving homeomorphisms} 

\subsection{Spaces of Radon measures} \mbox{} 

Suppose $Y$ is a locally connected, locally compact, $\sigma$-compact (separable metrizable)  space. 
Let ${\mathcal B}(Y)$ denote the $\sigma$-algebra of Borel subsets of $Y$.  A {\it Radon measure} on $Y$ is a measure $\mu$ on the measurable space $(Y, {\mathcal B}(Y))$ 
such that $\mu(K) < \infty$ for any compact subset $K$ of $Y$. 
Let ${\mathcal M}(Y)$ denote the set of Radon measures on $Y$. 
The {\it weak} topology $w$ on ${\mathcal M}(Y)$ is the weakest topology such
that the function 
$$\Phi_f : {\mathcal M}(Y) \lra {\Bbb R} \ : \ \mu \lmt \int_Y f \,d\mu$$
is continuous for any continuous function $f : Y \to {\Bbb R}$ with compact support. 
The set ${\mathcal M}(Y)$ is endowed with the weak topology $w$, otherwise specified. 

For $\mu \in {\mathcal M}(Y)$ and $A \in {\cal B}(Y)$, the restriction $\mu|_A$ is the Radon  measure on $A$ defined by $(\mu|_A)(B) = \mu(B)$ \ ($B \in {\mathcal B}(A)$). 

\begin{lemma}\label{lem_rest}{\rm (\cite[Lemma 2.2]{Be3})}
For any closed subset $A$ of $Y$, the map  
${\mathcal M}(Y) \to {\mathcal M}(A)  :  \mu \, \mapsto \, \mu|_A$  
is continuous at each $\mu \in {\mathcal M}(Y)$ with $\mu({\rm Fr}_MA) = 0$. 
\end{lemma}

We say that $\mu \in {\mathcal M}(Y)$ is {\it good} if 
$\mu(p) = 0$ for any point $p \in Y$ and $\mu(U) > 0$ for any nonempty open subset $U$ of $Y$.  
For $A \in {\mathcal B}(Y)$ let ${\mathcal M}_g^A(Y)$ denote the subspace of  ${\mathcal M}(Y)$ consisting of good Radon measures $\mu$ on $Y$ with $\mu(A) = 0$. 
For $\mu, \nu \in {\mathcal M}(Y)$, we say that 
$\nu$ is $\mu$-{\it biregular} if $\nu$ and $\mu$ have same null sets (i.e., $\nu(B)  = 0$ iff $\mu(B) = 0$ for any $B \in {\mathcal B}(Y)$). 
For $\mu \in {\mathcal M}_g^A(Y)$ we set
$${\mathcal M}_g^A(Y; \mu\reg) = 
\big\{ \nu \in {\mathcal M}_g^A(Y) \mid \nu \ \text{is} \ \mu\bireg \big\} \hspace{2mm} 
\text{(with the weak topology).}$$ 

For $h \in {\mathcal H}(Y)$ and $\mu \in {\mathcal M}(Y)$, the induced measures  $h_\ast \mu, h^\ast \mu \in {\mathcal M}(Y)$ are defined by 
\[ \mbox{$(h_\ast \mu)(B) = \mu(h^{-1}(B))$ \ \ \text{and} \ \ 
$(h^\ast \mu)(B) = \mu(h(B))$ \ \ 
($B \in {\mathcal B}(Y)$).} \] 
The group ${\mathcal H}(Y)$ acts  continuously on the space ${\mathcal M}(Y)$ by  
$h \cdot \mu = h_\ast \mu$. 
We say that $h \in {\mathcal H}(Y)$ is 
\begin{itemize}
\item[(i)\,] {\it $\mu$-preserving} if $h_\ast \mu = \mu$  (i.e., $\mu(h(B)) = \mu(B)$ for any $B \in {\mathcal B}(Y)$) and 
\item[(ii)] $\mu$-{\it biregular} if $h_\ast \mu$ and $\mu$ have the same null sets (i.e., $\mu(h(B)) = 0$ iff $\mu(B) = 0$ for any $B \in {\mathcal B}(Y)$). 
\end{itemize} 
Let ${\mathcal H}(Y; \mu) \subset {\mathcal H}(Y; \mu\reg)$ denote the subgroups of ${\mathcal H}(Y)$ consisting of $\mu$-preserving and $\mu$-biregular homeomorphisms of $Y$ respectively. 
For a subset $A$ of $Y$, 
the subgroups  
${\mathcal H}_A(Y; \mu)$, ${\mathcal H}_A(Y; \mu)_1$, ${\mathcal H}_{A, c}(Y; \mu)$, ${\mathcal H}_A(Y; \mu\reg)$, etc.  are defined according to the conventions in Sections 2.1 and 2.3. 

For spaces of embeddings, we use the following notations. 
Suppose $Y$ is a locally compact, $\sigma$-compact (separable metrizable)  space and $\mu \in {\cal M}(Y)$. 
For any $X \in {\cal B}(Y)$, an embedding $f : X \to Y$ is said to be  
\begin{itemize}
\item[(i)\ ] Borel if $f(X) \in {\cal B}(Y)$, 
\item[(ii)\,] $\mu$-biregular provided $f$ is Borel and $\mu(f(B)) = 0$ iff 
$\mu(B) = 0$ for any $B \in {\cal B}(X)$, 
\item[(iii)] $\mu$-preserving provided $f$ is Borel and $f : (X, \mu|_X) \cong (f(X), \mu|_{f(X)})$ is a measure preserving homeomorphism (i.e., $\mu(f(B)) = \mu(B)$ for any $B \in {\mathcal B}(X)$). 
\end{itemize} 
For a subset $A$ of $X$, let ${\cal E}_A(X, Y; \mu\reg)$ and  ${\cal E}_A(X, Y; \mu)$ denote the subspaces of 
${\cal E}_A(X, Y)$ consisting of $\mu$-biregular embeddings and $\mu$-preserving embeddings respectively. 

Suppose $M$ is a {\em compact connected} $n$-manifold and $\mu \in {\cal M}_g^\partial(M) \big( = {\cal M}_g^{\partial M}(M) \big)$. 

\begin{theorem}\label{thm_vNOU} {\rm (\cite{OU})} 
If $\nu \in {\mathcal M}_g^\partial(M)$ and $\nu(M) = \mu(M)$, then 
there exists $h \in {\mathcal H}_\partial(M)_1$ such that $h_\ast \mu = \nu$.
\end{theorem}

Let ${\cal M}_g^\partial(M; \mu) = 
\big\{ \nu \in {\cal M}_g^\partial(M; \mu\reg) \mid \nu(M) = \mu(M) \big\}$  
(with the weak topology). (See Section 3.2 for the definition 
in the case where $M$ is noncompact.) 
The group ${\cal H}(M; \mu\reg)$ acts continuously on 
${\cal M}_g^\partial(M; \mu)$ by $h \cdot \nu = h_\ast \nu$.  
This action induces the map 
$$\pi : {\cal H}(M; \mu\reg) \to {\cal M}_g^\partial(M; \mu) \ : \  h \lmt h_\ast \mu.$$ 

\begin{theorem}\label{thm_selection} $($\cite[Theorem 3.3]{Fa}$)$
The map $\pi$ admits a section 
$$\sigma : {\cal M}_g^\partial(M; \mu) \lra {\cal H}_\partial(M; \mu\reg)_1 \ \subset \ {\cal H}(M; \mu\reg) \hspace{2mm} \text{such that} \hspace{2mm} (\,\pi \sigma = id \hspace{2mm} \text{and}\,) \hspace{2mm} \sigma(\mu) = id_M.$$ 
\end{theorem} 

Next we recall basic facts on the product of measures. 
Suppose $(X, {\cal F}, \mu)$ and $(Y, {\cal G}, \nu)$ are 
$\sigma$-finite measure spaces. 
Let ${\cal F} \times {\cal G}$ denote  
the $\sigma$-algebra on $X \times Y$ generated by the family 
$\{ A \times B \mid A \in {\cal F}, B \in {\cal G} \}$. 
For $G \in {\cal F} \times {\cal G}$ and $x \in X$, the slice 
$G_x \subset Y$ is defined by 
$G_x = \{ y \in Y \mid (x,y) \in G \}$.  
It is well known that 
\begin{itemize}
\item[(1)] 
there exists a unique measure $\omega$ on the measurable space 
$(X \times Y, {\cal F} \times {\cal G})$ such that \\
\hspace{8mm} $\omega(A \times B) = \mu(A) \cdot \nu(B)$ \ \ 
$(A \in {\cal F}, B \in {\cal G})$ \ \ $($we follow the convention $0 \cdot \infty = 0$$)$, 
\item[(2)] for any $G \in {\cal F} \times {\cal G}$ \\
\hspace{5mm} 
(i) \ $\nu(G_x)$ $(x \in X)$ \ is an ${\cal F}$-measurable function on $X$ \ \ and \ \ 
(ii) \ $\ds \omega(G) = \int_X \nu(G_x)\,d\mu(x)$.  
\end{itemize} 
This result yields the following consequences on the product of Radon measures. 

\begin{proposition}\label{prop_prod_measure} 
Suppose $(X, \mu)$ and $(Y, \nu)$ are 
locally compact separable metrizable spaces with Radon measures. Then 
the following hold: 
\begin{itemize}
\item[(0)] ${\cal B}(X) \times {\cal B}(Y) = {\cal B}(X \times Y)$. 
\item[(1)] There exists a unique $\omega \in {\cal M}(X \times Y)$ such that 
$\omega(A \times B) = \mu(A) \cdot \nu(B)$ \ 
$(A \in {\cal B}(X), B \in {\cal B}(Y))$. 

\item[(2)] For any $G \in {\cal B}(X \times Y)$ 
\begin{itemize}
\item[(i)\,] $\nu(G_x)$ $(x \in X)$ \ is a ${\cal B}(X)$-measurable function on $X$ 
\ \ and \ \ {\rm (ii)} \ $\ds \omega(G) = \int_X \nu(G_x)\,d\mu(x)$. 
\end{itemize}
\end{itemize}
\end{proposition}

\noindent 
The measure $\omega$ is called the product of $\mu$ and $\nu$ and denoted by $\mu \times \nu$. 

\begin{proposition}\label{prop_prod_mp_homeo} 
Suppose $f : (X, \mu) \to (X_1, \mu_1)$ and $g : (Y, \nu) \to (Y_1, \nu_1)$ are 
homeomorphisms between 
locally compact separable metrizable spaces with Radon measures. Then the product homeomorphism \break 
$f \times g : (X \times Y, \mu \times \nu) \lra (X_1 \times Y_1, \mu_1 \times \nu_1)$ \ has the following properties: 
\begin{itemize}
\item[(1)] If $f$ and $g$ are biregular, then $f \times g$ is biregular. 
\item[(2)] If $f$ and $g$ are measure-preserving, then $f \times g$ is measure-preserving. 
\end{itemize}
\end{proposition}

\begin{proof}
For $G \in {\cal B}(X \times Y)$, we have \ \ (a) \ $\ds (\mu \times \nu)(G) = \int_X \nu(G_x)\,d\mu(x)$ \ \ and 
$$\hspace{5mm} 
{\rm (b)} \ \ \ds (\mu_1 \times \nu_1)\big( (f \times g)(G)\big) = 
\ds \int_{X_1} \nu_1\big(\big((f \times g)(G)\big)_{x_1}\big)\,d\mu_1(x_1) 
= \int_{X_1} \nu_1(g(G_{f^{-1}(x_1)}))\,d\mu_1(x_1).$$ 
\vskip 0mm 

(1) Note that 
\begin{itemize}
\item[(i)\,] 
$\ds (\mu \times \nu)(G) = 0$ \ iff 
\begin{tabular}[t]{l} 
$\nu(G_x) = 0$ \ ($\mu$-a.e. $x \in X$) \\[1mm]
\hspace{5mm} (i.e., \ $\exists$ $A \in {\cal B}(X)$ such that $\mu(A) = 0$ and $\nu(G_x) = 0$ ($x \in X - A$)), 
\end{tabular}
\vskip 1mm 
\item[(ii)] $\ds (\mu_1\times \nu_1)\big( (f \times g)(G)\big) = 0$ \ iff \ 
$\nu_1(g(G_{f^{-1}(x_1)})) = 0$ \ ($\mu_1$-a.e. $x_1 \in X_1$). 
\end{itemize} 
\vskip 1mm 
Since $f$ and $g$ are biregular, if (i) holds, then it follows that  
\[ \mbox{$f(A) \in {\cal B}(X_1)$, \ $\mu_1(f(A)) = 0$ \ and \ $\nu_1(g(G_{f^{-1}(x_1)})) = 0$ \ ($x_1 \in X_1 - f(A)$).} \] 
This implies (ii). The same argument shows the opposite implication. 
This means that $f \times g$ is biregular. 
\vskip 1mm 
(2) Since $f$ and $g$ are measure-preserving, it follows that 
$$\begin{array}[c]{ccl}
\ds (\mu_1 \times \nu_1)\big( (f \times g)(G)\big) &=& 
\ds \int_{X_1} \nu_1(g(G_{f^{-1}(x_1)}))\,d\mu_1(x_1) 
\ = \ \ds \int_{X_1} \nu(G_{f^{-1}(x_1)})\,d\mu_1(x_1) \\[4mm] 
&=& \ds \int_X \nu(G_x)\,d\mu(x) \ = \ (\mu \times \nu)(G).
\end{array}$$  
This means that $f \times g$ is measure-preserving. We also note that \ 
$(f \times g)^\ast(\mu_1 \times \nu_1) \in {\cal M}(X \times Y)$ \ satisfies the condition : for any $A \in {\cal B}(X)$ and $B \in {\cal B}(Y)$ 
$$\begin{array}[c]{ccl}
\big((f \times g)^\ast(\mu_1 \times \nu_1)\big)(A \times B) 
&=& (\mu_1 \times \nu_1)\big((f \times g)(A \times B)\big)
\ = \ (\mu_1 \times \nu_1)\big(f(A) \times g(B)\big) \\[2mm] 
&=& \mu_1(f(A)) \cdot \nu_1(g(B))
\ = \ \mu(A) \cdot \nu(B). 
\end{array}$$  
By definition we have 
$(f \times g)^\ast(\mu_1 \times \nu_1) = \mu \times \nu$. 
This also implies the conclusion.  
\end{proof} 

We conclude this subsection with some remarks on collars of the boundary of  
a submanifold. 
Suppose $M$ is an $n$-manifold and $\mu \in {\cal M}_g^\partial(M)$. 

\begin{remark}\label{rem_bdy}
Suppose $N$ is an $n$-submanifold of $M$ such that 
$\partial_+ N$ is compact. 
Since $\mu(\partial M) = 0$, we have $\mu(\partial N) = \mu(\partial_+ N)$. 
Take a bicollar $\partial_+ N \times [-1,1]$ of $\partial_+ N$ in $M$. 
Since $\partial_+ N \times [-1,1]$ is compact, it follows that 
$\mu(\partial_+ N \times [-1,1]) < \infty$ and 
$\big\{ t \in [-1,1] \mid \mu(\partial_+ N \times \{ t \}) \neq 0 \big\}$ is a countable subset of $[-1, 1]$. 
Hence, we can modify $N$ by adding or subtracting a thin collar of $\partial_+ N$ so that $\mu(\partial N) = \mu(\partial_+ N) = 0$.  
\end{remark}

Let $m$ denote the Lebesgue measure on the real line ${\Bbb R}$.  

\begin{lemma}\label{lem_collar} 
Suppose $N$ is an $n$-submanifold of $M$ such that $\partial_+ N$ is compact and $\mu(\partial_+ N) = 0$ and suppose $\nu \in {\cal M}_g^\partial(\partial_+ N)$. 
Then, there exists a bicollar $E = \partial_+ N \times [a, b]$ $(a < 0 < b)$ 
of $\partial_+ N$ in $M$  
such that $\partial_+ N = \partial_+ N \times \{ 0 \}$, $N \cap E = \partial_+ N \times [a, 0]$ 
and $\mu|_E = \nu \times (m|_{[a, b]})$. 
\end{lemma} 
  
\begin{proof}
Let ${\cal C}(\partial_+ N) = \{ F_1, \cdots, F_m \}$. 
For each $i = 1, \cdots, m$, 
choose a small bicollar $E_i = F_i \times [a_i, b_i]$ ($a_i < 0 < b_i$) such that 
$F_i = F_i \times \{ 0 \}$, $N \cap E_i = F_i \times [a_i, 0]$, 
$\mu(\partial_+ E_i) = 0$, $\mu(F_i \times [a_i, 0]) = |a_i|\nu(F_i)$ and 
$\mu(F_i \times [0, b_i]) = b_i\nu(F_i)$. 
We can apply Theorem~\ref{thm_vNOU} to $$\mu|_{F_i \times [a_i, 0]}, 
\nu|_{F_i} \times (m|_{[a_i, 0]}) \in {\cal M}_g^\partial(F_i \times [a_i, 0]) \ \ \text{ and } \ \ 
\mu|_{F_i \times [0, b_i]}, 
\nu|_{F_i} \times (m|_{[0, b_i]}) \in {\cal M}_g^\partial(F_i \times [0, b_i])$$
to replace the identification of the collar $E_i = F_i \times [a_i, b_i]$ so that 
$\mu|_{E_i} = \nu|_{F_i} \times (m|_{[a_i, b_i]})$. 
Finally, take $a, b$ such that $\max_i a_i < a < 0 < b < \min_i b_i$ and 
set $E = \partial_+ N \times [a,b] = \bigcup_i \big( F_i \times [a,b] \big)$. 
\end{proof} 


\subsection{End compactification and finite-end weak topology} (cf.\,\cite{AP2, Be3}) 
 
In order to extend the selection theorem~\ref{thm_selection} to the noncompact case, 
it is necessary to include the information of the ends.  
Suppose $Y$ is a noncompact, connected, locally connected, locally compact, separable metrizable space. 
Let ${\cal K}(Y)$ denote the collection of all compact subsets of $Y$. 
An {\em end} of $Y$ is a function $e$ which assigns an $e(K) \in {\cal C}(Y - K)$ to each $K \in {\cal K}(Y)$ such that 
$e(K_1) \supset e(K_2)$ if $K_1 \subset K_2$. 
The set of ends of $Y$ is denoted by $E_Y$. 
The {\em end compactification} of $Y$ is the space $\overline{Y} = Y \cup E_Y$ 
equipped with the topology defined by the following conditions: 
(i) $Y$ is an open subspace of $\overline{Y}$, 
(ii) the fundamental open neighborhoods of $e \in E_Y$ are given by 
\[ N(e, K) = e(K) \,\cup \,\{ e' \in E_Y \mid e'(K) = e(K)\} \hspace{4mm} (K \in {\cal K}(Y)). \] 
Then $\overline{Y}$ is a connected, locally connected, compact, metrizable space, 
$Y$ is a dense open subset of $\overline{Y}$ and  
$E_Y$ is a compact 0-dimensional subset of $\overline{Y}$. 

For $h \in {\cal H}(Y)$ and $e \in E_Y$ we define $h(e) \in E_Y$ by $h(e)(K) = h(e(h^{-1}(K)))$ $(K \in {\cal K}(Y))$. 
Each $h \in {\cal H}(Y)$ has a unique extension $\overline{h} \in {\cal H}(\overline{Y})$ defined by 
$\overline{h}(e) = h(e)$ $(e \in E_Y)$. 
The map \ ${\cal H}(Y) \to {\cal H}(\overline{Y})$ : $h \mapsto \overline{h}$ \ is a continuous group homomorphism.  
For $A \subset Y$ we set ${\cal H}_{A \cup E_Y}(Y) = \{ h \in {\cal H}_A(Y) \mid \overline{h}|_{E_Y} = id_{E_Y}\}$. 
Note that ${\cal H}_{A \cup E_Y}(Y)_0 = {\cal H}_{A}(Y)_0$. 

Let $\mu \in {\mathcal M}(Y)$. An end $e \in E_Y$ is said to be $\mu$-{\it finite} if 
$\mu(e(K)) < \infty$ \ for some $K \in {\mathcal K}(Y)$. 
Let $E_Y^\mu = \{ e \in E_Y \mid e \text{ is } \mu\text{-finite}\}$. 
Then $Y \cup E_Y^\mu$ is an open subset of $\overline{Y}$.  
For $A \in {\mathcal B}(Y)$ and $\mu \in {\mathcal M}_g^A(Y)$ we set 
$$\begin{array}[c]{l}
{\mathcal M}_g^A(Y; \mu\ereg) 
= \big\{ \nu \in {\mathcal M}_g^A(Y) \mid 
\nu \text{ is } \mu\bireg, \ E_Y^\nu = E_Y^\mu \big\}, \\[2.5mm]
{\mathcal M}_g^A(Y; \mu) 
= \big\{ \nu \in {\mathcal M}_g^A(Y; \mu\ereg) \mid 
\nu(Y) = \mu(Y)\big\}.
\end{array}$$
The {\it finite-ends weak} topology $ew$ on ${\mathcal M}_g^A(Y; \mu\ereg)$ is 
the weakest topology such that the function 
$$\Phi_f : {\mathcal M}_g^A(Y; \mu\ereg) \lra {\Bbb R} \ : \ \nu \lmt \int_Y f|_Y \,d\nu$$
is continuous for any continuous function $f : Y \cup E_Y^\mu \to {\Bbb R}$ with compact support. 

There is an alternative description of this topology (\cite[\S3, p\,245]{Be3}). 
Consider the space ${\cal M}(Y \cup E_Y^\mu)$ (with the weak topology). 
Each $\nu \in {\mathcal M}_g(Y; \mu\ereg)$ has a natural extension $\overline{\nu} \in {\mathcal M}_g(Y \cup E_Y^\mu)$ 
defined by $\overline{\nu}(B) = \nu(B \cap Y)$ ($B \in {\mathcal B}(Y \cup E_Y^\mu)$). 
The topology $ew$ on ${\mathcal M}_g^A(Y; \mu\ereg)$ is 
the weakest topology for which the injection 
$$\iota : {\mathcal M}_g^A(Y; \mu\ereg) \lra {\mathcal M}(Y \cup E_Y^\mu)_w  \ : \ \nu \, \longmapsto \, \overline{\nu}$$ 
is continuous.
The symbol ${\mathcal M}_g^A(Y; \mu\ereg)_{ew}$ denotes the space ${\mathcal M}_g^A(Y; \mu\ereg)$ endowed with the topology $ew$. 

We say that $h \in {\mathcal H}(Y)$ is $\mu$-{\it end-biregular} if $h$ is $\mu$-biregular and $E_Y^{h_\ast \mu} = E_Y^\mu$ (i.e., $\overline{h}(E_Y^\mu) = E_Y^\mu$). 
Let ${\mathcal H}(Y; \mu\ereg)$ denote the subgroup of ${\mathcal H}(Y)$ consisting of $\mu$-end-biregular  homeomorphisms of $Y$. 

Suppose $M$ is a connected $n$-manifold and $\mu \in {\mathcal M}_g^\partial(M)$. 
The group ${\mathcal H}(M; \mu\ereg)$ acts continuously on ${\mathcal M}_g^\partial(M; \mu)_{ew}$ 
by $h \cdot \nu = h_\ast \nu$. 
This action induces the map 
$$\pi : {\mathcal H}(M; \mu\ereg) \lra {\mathcal M}_g^\partial(M; \mu)_{ew} \ : \  h \lmt h_\ast \mu.$$ 

\begin{theorem}\label{thm_selection_noncpt} {\rm (\cite[Theorem 4.1]{Be3})} The map $\pi$ has a section 
$$\sigma : {\mathcal M}_g^\partial(M; \mu)_{ew} \lra {\cal H}_\partial(M; \mu\ereg)_1 \ = \ {\mathcal H}_\partial(M; \mu\reg)_1 
\hspace{2mm} \text{such that} \hspace{2mm} (\,\pi \sigma = id \hspace{2mm} \text{and}\,) \hspace{2mm} \sigma(\mu) = id_M.$$ 
\end{theorem}

\subsection{End charge homomorphism} \mbox{} 

We recall basic properties of the end charge homomorphisms defined in \cite[Section 14]{AP2}. 
Suppose $Y$ is a connected, locally connected, locally compact separable, metrizable space. 
Let ${\cal Q}(E_Y)$ denote the algebra of clopen subsets of $E_Y$ and 
let ${\cal B}_c(Y) = \{ C \in {\cal B}(Y) \mid {\rm Fr}_Y\,C \mbox{ is compact} \}$.  
For each $C \in {\cal B}_c(Y)$ let 
\[ E_C = \{ e \in E_Y \mid e(K) \subset C \mbox{ for some } K \in {\cal K}(Y)\}
\hspace{5mm} \text{and} \hspace{5mm} \overline{C} = C \cup E_C \subset \overline{Y}. \]  
Note that 
(i) $E_C \in {\cal Q}(E_Y)$ and $\overline{C}$ is a neighborhood of $E_C$ in $\overline{Y}$ with $\overline{C} \cap E_Y = E_C$, 
(ii) for $C, D \in {\cal B}_c(Y)$ it follows that $E_C = E_D$ iff 
$C \Delta D = (C - D) \cup (D - C)$ is relatively compact (i.e., has the compact closure) in $Y$, 
(iii) if $C \in {\cal B}_c(Y)$ and $h \in {\cal H}_{E_Y}(Y)$, 
then $h(C) \in {\cal B}_c(Y)$ and $E_{h(C)} = E_C$. 

An {\em end charge} of $Y$ is a finitely additive signed measure $c$ on ${\cal Q}(E_Y)$, that is, 
a function $c : {\cal Q}(E_Y) \to {\Bbb R}$ which satisfies the following condition:  
\[ \mbox{$c(F \cup G) = c(F) + c(G)$ \ for \ $F, G \in {\cal Q}(E_Y)$ \ with \ $F \cap G = \emptyset$.} \]  
Let ${\cal S}(Y)$ denote the space of end 
charges $c$ of $Y$ endowed with the {\em weak topology} (or the product topology). 
This topology is the weakest topology such that the function 
\[ \Psi_F : {\cal S}(Y) \lra {\Bbb R} \ : \ c \longmapsto c(F) \] 
is continuous for any $F \in {\cal Q}(E_Y)$. 
For $\mu \in {\cal M}(Y)$ let  
\[ 
\mbox{${\cal S}(Y, \mu) = \big\{ c \in {\cal S}(Y) \mid$ (i) $c(F) = 0$ for 
$F \in {\cal Q}(E_Y)$ with $F \subset E_Y^\mu$ \ and \ (ii) $c(E_Y) = 0 \,\big\}$} \] 
(with the weak topology). 
Then ${\cal S}(Y)$ is a topological linear space and ${\cal S}(Y, \mu)$ is a linear subspace. 

For $h \in {\cal H}_{E_Y}(Y; \mu)$ 
the end charge $c_h^\mu \in {\cal S}(Y, \mu)$ is defined as follows: 
For any $F \in {\cal Q}(E_Y)$ there exists $C \in {\cal B}_c(Y)$ with $E_C = F$.  
Since $\overline{h}|_{E_Y} = id$, it follows that $E_C = E_{h(C)}$ and that 
$C \Delta \,h(C)$ is relatively compact in $Y$. 
Thus $\mu(C - h(C)), \mu(h(C) - C) < \infty$ and we can define  
\[ c_h^\mu(F) = \mu(C - h(C)) - \mu(h(C) - C) \in {\Bbb R}.  \] 
This quantity is independent of the choice of $C$. 

\begin{proposition} The end charge homomorphism $c^\mu : {\cal H}_{E_Y}(Y; \mu) \lra {\cal S}(Y, \mu)$ is   
a continuous group homomorphism $($\cite[Section 14.9, Lemma 14.21\,(iv)]{AP2}$)$. 
\end{proposition} 

In \cite{Ya2} we have shown that, 
for any connected $n$-manifold $M$ and $\mu \in {\cal M}_g^\partial(M)$, 
the end charge homomorphism $c^\mu : {\cal H}_{E_M}(M; \mu) \to {\cal S}(M; \mu)$ has a (non-homomorphic) section 
$s :  {\cal S}(M, \mu) \to {\cal H}_\partial(M; \mu)_1$. 

For any subset $A$ of $Y$ we have the restriction of $c^\mu$ 
$$c^\mu_A : {\cal H}_{A \cup E_Y}(Y; \mu) \to S(Y, \mu).$$

The kernel of the homomorphism $c^\mu$ is denoted by ${\rm ker}\,c^\mu$. 
Note that ${\cal H}_c(M; \mu) \subset {\rm ker}\,c^\mu$ and 
$\big({\rm ker}\,c^\mu\big)_A = {\rm ker}\,c^\mu_A$. 
By the definition, if $h \in {\rm ker}\,c^\mu$, then for any $C \in {\cal B}_c(Y)$ we have 
$\mu(C - h(C)) = \mu(h(C) - C)$. 

\begin{lemma}\label{lem_ker} Suppose $h \in {\rm ker}\,c^\mu$ and $C \in {\cal B}_c(Y)$. 
If $L \in {\cal B}(C \cap h(C))$ and $C - L$ is relatively compact in $Y$, then 
$h(C) - L$ is also relatively compact and 
$\mu(h(C) - L) = \mu(C - L)$. 
\end{lemma}

\begin{proof} 
Since $\mu(C - h(C)) = \mu(h(C) - C)$, the assertion follows from the equalities : 
$$h(C) - L = (h(C) - C) \cup \big((C \cap h(C)) - L\big) \hspace{4mm} \text{and} \hspace{4mm}
C - L = (C - h(C)) \cup \big((C \cap h(C)) - L\big).$$  
\vskip -6mm 
\end{proof} 


\section{Weak extension theorem for biregular homeomorphisms} 

Throughout this section $M$ is an $n$-manifold and 
$\mu \in {\cal M}_g^{\partial}(M)$. 
The weak extension theorem for the group 
$G = {\cal H}(M; \mu\reg)$ is already obtained in 
\cite{Fa}. 
In this section we discuss some consequences of this extension theorem. 
In Section 5 we combine 
the weak extension theorem for ${\cal H}(M; \mu\reg)$ and 
the selection theorem for $\mu$-biregular measures (Theorems~\ref{thm_selection} and \ref{thm_selection_noncpt}) in order to obtain 
the weak extension theorems for the groups ${\cal H}(M; \mu)$ and ${\rm ker}\,c^\mu$. 

First we recall the deformation theorem for $\mu$-biregular embeddings \cite[Theorem 4.1]{Fa}. 
For $X \in {\cal B}(M)$ and $A \subset X$, 
let ${\cal E}^\ast_A(X, M; \mu\reg)$ denote the space of proper $\mu$-biregular embeddings $f : X \to M$ with $f|_A = id_A$, endowed with the compact-open topology (cf. Sections 2.1 and 3.1). 

Suppose 
$C$ is a  compact subset of $M$, 
$U \in {\cal B}(M)$ is a neighborhood of $C$ in $M$ and 
$D \subset E$ are two closed subsets of $M$ such that $D \subset {\rm Int}_M E$. 

\begin{theorem}\label{thm-deform_bireg} $($\cite[Theorem 4.1]{Fa}$)$  
For any compact neighborhood $K$ of $C$ in $U$,  
there exists a neighborhood $\U$ of $i_U$ in ${\mathcal E}_{E \cap U}^\ast(U, M; \mu\reg)$ 
and a homotopy \ 
$\varphi : \U \times [0, 1] \lra {\mathcal E}^\ast_{D \cap U}(U, M; \mu\reg)$ \  such that 
\begin{itemize}
\item[{\rm (1)}]  for each $f \in {\cal U}$, \\[1mm]
\hspace{5mm} {\rm (i)} \ \ $\varphi_0(f) = f$, \hspace{3mm} {\rm (ii)} 
$\varphi_1(f)|_C = i_C$, \hspace{3mm} {\rm (iii)} 
$\varphi_t(f)|_{U - K} = f|_{U - K}$ $(t \in [0,1])$, \\[0.8mm] 
\hspace{4mm} {\rm (iv)} \ if $f = id$ on $\partial_- U$, then $\varphi_t(f) = id$ on $\partial_- U$ $(t \in [0,1])$, 
\vskip 1mm 
\item[{\rm (2)}] $\varphi_t(i_U) = i_U$ $(t \in [0,1])$. 
\end{itemize}
\end{theorem} 

Theorem~\ref{thm-deform_bireg} is equivalent to the next weak extension  theorem. 

\begin{theorem} \label{thm_ext_bireg}$($\cite[Corollary 4.2]{Fa}$)$ 
For any compact neighborhood $L$ of $C$ in $U$,  
there exists a neighborhood ${\cal U}$ of $i_U$ in  ${\cal E}^\ast_{E \cap U}(U, M; \mu\reg)$ and a homotopy \ 
$s : {\cal U} \times [0,1] \to {\cal H}_{D \cup (M - L)}(M; \mu\reg)_1$ \ such that 
\begin{itemize}
\item[(1)] for each $f \in {\cal U}$ \\
\hspace{3mm} {\rm (i)} $s_0(f) = id_M$, \hspace{3mm} {\rm (ii)} 
$s_1(f)|_C = f|_C$, \hspace{3mm} {\rm (iii)} 
if $f = id$ on $\partial_- U$, then $s_t(f) = id$ on $\partial M$, 
\item[{\rm (2)}] $s_t(i_U) = id_M$ $(t \in [0,1])$. 
\end{itemize} 
\end{theorem} 

\noindent (In \cite[Corollary 4.2]{Fa} the map $s_1$ alone is mentioned.)  

Now we discuss some consequences of Theorem \ref{thm_ext_bireg} 
for the group $G = {\cal H}(M; \mu\reg)$. 
Suppose $X$ is a compact subset of $M$. Note that 
$G_X = {\cal H}_X(M; \mu\reg)$. 

Suppose $C$ is a compact subset of $M$ with $X \subset C$ 
and $U$ is a neighborhood of $C$ in $M$. 
Consider the pull-back diagram : \\[-2mm] 
\hspace{30mm} 
$\begin{array}[c]{cccl}
& p' & & \\
p^\ast G_X & \lra & G_X \\[0.5mm] 
\pi' \left\downarrow \makebox(0,10){} \right. \ \ & & \ \ \left\downarrow \makebox(0,10){}\right. \pi & \hspace{3mm} \text{,where} \ \ \ \pi(h) = h|_C \ \ \ \text{and} \ \ \ p(f) = f|_C. \\[2mm] 
{\cal E}^G_X(U, M) & \lra & {\cal E}^G_X(C, M) & \\ 
& p & & 
\end{array}$
\vskip 2mm 
By Theorem \ref{thm_ext_bireg} the pair $(U, C)$ has WEP$_G$. 
Hence it has LSP$_G$ and also LSP$_{G_X}$. 
Thus the next assertion follows from Lemma~\ref{lem_p-bdl_2}. 

\begin{lemma}\label{lem_p-bdl_bireg} 
The induced map $\pi' : p^\ast G_X \to {\cal E}^G_X(U, M)$ is a principal $G_C$-bundle. 
\end{lemma} 

Suppose $N$ is a compact $n$-submanifold of $M$ such that 
$\mu(\partial_+ N) = 0$ and $X \subset {\rm Int}_M N$. 
Take any compact $n$-submanifold $N_1$ of $M$ such that 
$\mu(\partial_+ N_1) = 0$ and 
$N_1$ is obtained from $N$ by adding an outer collar of $\partial_+ N$. 
We obtain the pull-back diagram: 
$$\begin{array}[c]{cccl}
& p' & & \\
p^\ast G_X & \lra & G_X & \\[0.5mm] 
\pi' \left\downarrow \makebox(0,10){} \right. \ \ & & \ \ \left\downarrow \makebox(0,10){}\right. \pi & \hspace{8mm} 
\smash{\raisebox{2mm}{,where \ \ \ $\pi(g) = g|_N, \ \ \ p(f) = f|_N$ \ \ \ and}} \\[2mm] 
{\cal E}^G_X(N_1, M) & \lra & {\cal E}^G_X(N, M) & \hspace{22.5mm} \, p^{-1}(i_N) = {\cal E}_N^G(N_1, M). \\ 
& p & & 
\end{array}$$ 

\begin{lemma}\label{lem_path_bireg} \mbox{} 
There exists a path $h : [0,1] \to G_X$ such that 
\[ \mbox{$h_0 = id_M$, \ \ $h_1(N_1) = N$ \ \ and \ \ $h_t(N_1) \subset N_1$, \ $h_t(N) \subset N$ \ $(t \in [0,1])$.} \]  
\end{lemma} 

\begin{proof} 
(1) Let $m$ denote the Lebesgue measure on ${\Bbb R}$.
We can find a bicollar 
$E = \partial_+ N \times [a, b]$ ($a < 0$, $b > 1$) of $\partial_+ N$ in $M - X$ 
and $\nu \in {\cal M}_g^\partial(\partial_+ N)$ such that 
\[ \mbox{(i) $\partial_+ N = \partial_+ N \times \{ 0 \}$, \ 
$\partial_+ N_1 = \partial_+ N \times \{ 1 \}$ \ \ and \ \ (ii) 
$\mu|_E = \nu \times (m|_{[a,b]})$.} \]  
This follows from the following observation. 
First take any bicollar $E' = \partial_+ N \times [-1, 2]$ of $\partial_+ N$ in $M - X$ which satisfies (i) and the weaker condition (ii)$'$ $\mu(\partial_+ N \times \{ -1 \}) = \mu(\partial_+ N \times \{ 2 \}) = 0$. 
Let ${\cal C}(\partial_+ N) = \{ F_1, \dots, F_m \}$ and 
set $E_i' = F_i \times [-1, 2]$ ($i = 1, \cdots, m$). 
Choose any $\nu \in {\cal M}_g^\partial(\partial_+ N)$ such that 
$\nu(F_i) = \mu(F_i \times [0,1])$ ($i = 1, \cdots, m$).  
For each $i = 1, \cdots, m$, determine $a_i < 0$ and $b_i > 1$ by 
$|a_i| \nu(F_i) = \mu(F_i \times [-1, 0])$ and 
$(b_i - 1) \nu(F_i) = \mu(F_i \times [1, 2])$, and  
reparametrize $F_i \times [-1,0]$ to $F_i \times [a_i,0]$ and 
$F_i \times [1,2]$ to $F_i \times [1, b_i]$. 
We can apply Theorem~\ref{thm_vNOU} on $F_i \times [a_i,0]$, 
$F_i \times [0,1]$ and $F_i \times [1,b_i]$ to obtain 
a new identification $E_i' = F_i \times [a_i, b_i]$ so that 
$\mu|_{E_i'} = \nu \times (m|_{[a_i, b_i]}).$ 
Take $a, b$ such that $\max_i a_i < a < 0$ and $1 < b < \min_i b_i$, 
and set $E = \bigcup_i (F_i \times [a,b])$. 

(2) Choose $\lambda \in {\cal H}_\partial([a, b])$ such that 
$\lambda$ is piecewise affine and $\lambda(0) = a/2$, $\lambda(1) = 0$. 
We obtain two isotopies \\[-2mm] 
\hspace{17mm} 
\begin{tabular}[t]{llll}
$\lambda_t \in {\cal H}_\partial([a, b])$ & ($t \in [0,1]$) & defined by & $\lambda_t(s) = (1-t)s + t \lambda(s)$ \ \ and \\[1mm]
$g_t \in {\cal H}_{ \partial_+ N \times \{ a, b \}}(\partial_+ N \times [a, b])$ & ($t \in [0,1]$) & defined by & $g_t(y, s) = (y, \lambda_t(s))$. 
\end{tabular} \\[2.5mm] 
Note that $\lambda_0 = id$, $\lambda_1([a,1]) = [a, 0]$, $\lambda_t([a,0]) \subset [a,0]$ and $\lambda_t([a,1]) \subset [a,1]$. 
Since $\lambda_t$ is also piecewise affine, it is seen that $\lambda_t$ is $m|_{[a,b]}$-biregular. Then 
each $g_t$ is $\nu \times (m|_{[a,b]})$-biregular by Proposition~\ref{prop_prod_mp_homeo}.  
Finally, the required isotopy \ $h_t \in {\cal H}_{E^c}(M; \mu\reg) \subset G_X$ \ ($t \in [0,1]$) \ is defined by $h_t|_E = g_t$. 
\end{proof}

By Lemmas~\ref{lem_p-bdl_bireg}, \ref{lem_path_bireg} and \ref{lem_he} we have the following conclusions. 

\begin{lemma} 
{\rm (1)} The induced map $\pi' : p^\ast G_X \to {\cal E}^G_X(N_1, M)$ is a principal $G_N$-bundle. 
\begin{itemize}
\item[(2)] The map $p : {\cal E}^G_X(N_1, M) \to {\cal E}^G_X(N, M)$ is a homotopy equivalence. 

\item[(3)] There exists a strong deformation retraction $\chi_t$ $(t \in [0,1])$ of 
${\cal E}_N^G(N_1, M)$ onto the singleton $\{ i_{N_1} \}$. 

\item[(4)]The map $p' :  p^\ast G_X \to G_X$ is a homotopy equivalence. 
\end{itemize}
\end{lemma}

\begin{corollary}\label{cor_fib_bireg} 
Suppose $X$ is a compact subset of $M$ and 
$N$ is a compact $n$-submanifold of $M$ such that 
$\mu(\partial N) = 0$ and $X \subset {\rm Int}_M N$. 
Then the restriction map \vspace{-1mm} 
$${\cal H}_N(M; \mu\reg) \ \subset \ {\cal H}_X(M; \mu\reg) \lraw[\pi] {\cal E}^{{\cal H}(M; \mu\reg)}_X(N, M) \hspace{8mm} \text{defined by \hspace{3mm} $\pi(h) = h|_{N}$}$$   
is a fibration up to homotopy equivalences and has the exact sequence for  homotopy groups. 
\end{corollary} 


\section{Weak extension theorem for measure-preserving homeomorphisms} 

Throughout this section  
$M$ is an $n$-manifold and $\mu \in {\cal M}_g^{\partial}(M)$. 
In this section we combine 
the weak extension theorem for $G = {\cal H}(M; \mu\reg)$ (Theorem~\ref{thm_ext_bireg}) and 
the selection theorem for $\mu$-biregular measures (Theorems~\ref{thm_selection} and \ref{thm_selection_noncpt}) in order to obtain  
the weak extension theorems for the groups $H = {\cal H}(M; \mu)$ and $F = {\rm ker}\,c^\mu$. 
We also discuss 
a non-ambient weak deformation of measure-preserving embeddings (Theorem~\ref{thm_non-amb_def}). 
Some application to the group $H_c = {\cal H}_c(M; \mu)$ endowed with the Whitney  topology is provided in Section 6. 

\subsection{Weak extension theorem for ${\cal H}(M; \mu)$} \mbox{} 

We obtain the weak extension theorem for ${\cal H}(M; \mu)$ in a general form (Theorem 5.1, cf. \cite[Theorem 4.12]{Fa}). 
This answers Problem~\ref{prob_large} and also 
leads us to the weak extension theorem for ${\rm ker}\,c^\mu$ 
in Section 5.2. 
(Recall that $M$ is an $n$-manifold, $\mu \in {\cal M}_g^{\partial}(M)$, $G = {\cal H}(M; \mu\reg)$ and $H = {\cal H}(M; \mu)$.) 

For $A, B \in {\cal B}(M)$, consider 
the subset $G^{A, B}$ of $G$ defined by 
$$G^{A, B} = \big\{ h \in G \mid h|_A \in {\cal E}(A, M; \mu) 
\ \text{and} \ \mu(h(L)) = \mu(L) \ (L \in {\cal C}(M - B)) \big\}.$$ 
When $A = B$, we simply write $G^{A}$. 
For any $X \subset M$ we have the pair $(G^{A, B}_X, H_X)$ of subsets in $G_X$. 

\begin{lemma}\label{lem_ext_mp}
Suppose $N$ is a compact $n$-submanifold of $M$ with $\mu(\partial N) = 0$, $U \in {\cal B}(M)$ is a neighborhood of $N$ in $M$ and $X$ is a closed subset of $\partial M$ with $X \cap N = \emptyset$. 
Then the triple $(M, U, N)$ has WEP$(G_X^N, H_X)$. 
\end{lemma} 

\begin{proof} 
{Case 1:} First we consider the case where $M$ is connected. 

Since ${\cal E}^{G_X^N}(U, M) \subset {\cal E}^\ast(U, M; \mu\reg)$, 
by Theorem \ref{thm_ext_bireg} applied to $(U, C) = (M - X, N)$, 
there exists a neighborhood $\U$ of $i_U$ in ${\cal E}^{G_X^N}(U, M)$ and a map 
$\sigma : \U \times [0,1] \to (G_X)_1$ 
such that 
\begin{itemize}
\item[(i)\,] for each $f \in {\cal U}$ \\
\hspace{3mm} {\rm (a)} $\sigma_0(f) = id_M$, \hspace{3mm} {\rm (b)} 
$\sigma_1(f)|_N = f|_N$, \hspace{3mm} {\rm (c)} 
if $f = id$ on $\partial_- U$, then $\sigma_t(f) = id$ on $\partial M$, 
\vskip 0.5mm 
\item[(ii)] $\sigma_t(i_U) = id_M$ $(t \in [0,1])$. 
\end{itemize} 

(1) First we modify the map $\sigma$ to achieve the following additional condition: \hspace{3mm} (i)\,(b$'$) \ $\sigma_1(f) \in H$. \\[1mm] 
Consider the induced map \hspace{5mm} 
$\nu : \U \times [0,1] \lra {\cal M}_g^{\partial}(M; \mu)_{ew} \hspace{5mm} \text{defined by} \hspace{5mm} \nu_t(f) = \sigma_t(f)^\ast \mu.$ \\[1mm]
Since $M$ is connected, each $L \in {\cal C}(N^c)$ meets $\partial_+ N$. 
Since $\partial_+ N$ is compact, it follows that ${\cal C}(N^c)$ is a finite set. 
We note that $\nu_1(f)|_{L} \in {\cal M}_g^{\partial}(L; \mu|_{L})$
for any $f \in {\cal U}$ and $L \in {\cal C}(N^c)$. 
In fact, since $\nu_1(f) \in {\cal M}_g^{\partial}(M; \mu\ereg)$ and $\mu(\partial N) = 0$,  
we have $\nu_1(f)|_{L} \in {\cal M}_g^{\partial}(L; \mu|_L\ereg)$.
It remains to show that $\nu_1(f)(L) = \mu(L)$. 
Since $f \in {\cal E}^{G_X^N}(U, M)$, 
there exists $h \in G_X^N$ such that $f = h|_U$. 
Then $k \equiv h^{-1}\sigma_1(f) \in {\cal H}_N(M)$. 
Since $M$ is connected, we see that   
$N \cap L \neq \emptyset$, and since 
$k = id$ on $N$, we have $k(L) = L$.  
Hence, $\sigma_1(f)(L) = h(L)$ and  
it follows that \hspace{5mm} $\nu_1(f)(L) = \mu(\sigma_1(f)(L)) = \mu(h(L )) = \mu(L)$.
\vskip 1mm 
For each $L \in {\cal C}(N^c)$ we obtain the map \hspace{5mm} 
$\U \lra {\cal M}_g^{\partial}(L; \mu|_{L})_{ew} \ : \ f \longmapsto \nu_1(f)|_{L}.$ \\[1mm] 
By the alternative description of the finite-ends weak topology and Lemma~\ref{lem_rest}, this map is seen to be continuous (cf. \cite[Lemma 3.2]{Ya1}).  
By Theorem \ref{thm_selection_noncpt} there exists a map 
\[ \mbox{$\eta_L : {\cal M}_g^\partial(L; \mu|_{L})_{ew} \lra 
{\cal H}_\partial(L;\mu|_{L}\reg)_1$ \hspace{2mm} 
such that \hspace{2mm} 
$\eta_L(\nu)_\ast(\mu|_{L}) = \nu$ \ \ and \ \ $\eta_L(\mu|_{L}) = id_{L}$.} \]  
Define the map \hspace{5mm} 
$\tau_L : \U \times [0,1] \to {\cal H}_{\partial}(L; \mu|_{L}\reg)_1$ \hspace{3mm} by \hspace{3mm} 
$\tau_L(f, t) = \eta_L((1-t)\mu|_{L} + t \nu_1(f)|_{L}).$
\vskip 1mm 
Combining $\tau_L$ ($L \in {\cal C}(N^c)$), we obtain the map 
$$\tau : \U \times [0,1] \to {\cal H}_{N \cup \partial M}(M; \mu\reg)_1 \hspace{4mm} \text{defined by} \hspace{4mm} 
\tau(f,t) = 
\left\{ \hspace{-1mm} 
\begin{array}[c]{ll}
\tau_L(f, t)& \text{on $L \in {\cal C}(N^c)$} \\[2mm]
\ \ id & \text{on $N$}.  
\end{array} \right.$$
Note that $\tau_0(f) = id_M$ and $\tau_1(f)_\ast \mu = \nu_1(f)$. 
Define a map 
$$\sigma' : \U \times [0,1] \lra {\cal H}_X(M; \mu\reg)_1 \hspace{5mm} \text{by} \hspace{5mm}
\sigma'_t(f) = 
\left\{ \hspace{-1mm} 
\begin{array}[c]{ll}
\sigma_{2t}(f) & (t \in [0, 1/2]) \\[2mm]
\sigma_1(f) \tau_{2t-1}(f) & (t \in [1/2,1]).  
\end{array} \right.$$
\vskip 1mm 
\noindent Then the map $\sigma'$ satisfies the conditions (i)\,(a), (b), (c) and (ii). 
The condition (i)\,(b$'$) is verified by 
$$\sigma'_1(f)_\ast \mu \, = \, \sigma_1(f)_\ast \tau_1(f)_\ast \mu 
\, = \, \sigma_1(f)_\ast \nu_1(f) \, = \, \sigma_1(f)_\ast \sigma_1(f)^\ast \mu \, = \, \mu.$$ 

(2) To see that the triple $(M, U, N)$ has WEP$(G^N_X, H_X)$, we construct a map 
$s : {\cal U} \times [0,1] \to H_X$ \break 
such that  
\begin{itemize}
\item[(iii)] for each $f \in {\cal U}$ \\
\hspace{3mm} {\rm (a)} $s_0(f) = id_M$, \hspace{3mm} {\rm (b)} 
$s_1(f)|_N = f|_N$, \hspace{3mm} {\rm (c)} 
if $f = id$ on $\partial_- U$, then $s_t(f) = id$ on $\partial M$, 
\vskip 0.5mm 
\item[(iv)] $s_t(i_U) = id_M$ $(t \in [0,1])$. 
\end{itemize} 
Consider the induced map \hspace{5mm} $\nu' : \U \times [0,1] \lra {\cal M}_g^\partial(M; \mu)_{ew} \hspace{4mm}  \text{defined by} \hspace{4mm} \nu'_t(f) = \sigma'_t(f)^\ast \mu$. \\[1mm]  
It is seen that $\nu'_0(f) = \nu'_1(f) = \mu$. 
By Theorem \ref{thm_selection_noncpt} there exists a map  
\[ \mbox{$\eta : {\cal M}_g^\partial(M; \mu)_{ew} \to (G_\partial)_1$ \hspace{2mm}  
such that \hspace{2mm} $\eta(\nu)_\ast \mu = \nu$ \ \ and \ \ 
$\eta(\mu) = id_M$.} \] 
The required map $s$ is defined by \hspace{10mm} 
$s_t(f) = \sigma'_t(f) \eta(\nu'_t(f)) \ \ ((f,t) \in {\cal U} \times [0,1]).$ \\[1mm] 
The conditions (iii) and (iv) are easily verified. For example, 
(iii) (b) is seen by 
\[ s_1(f) = \sigma'_1(f) \eta(\nu'_1(f)) = \sigma'_1(f) \eta(\mu) = \sigma'_1(f) \hspace{4mm} \text{and} \hspace{4mm} s_1(f)|_N = \sigma'_1(f)|_N = f|_N. \] 

\noindent {Case 2:} Next we treat the general case where $M$ may not be connected. 

By Lemma~\ref{l_WEP} we may assume that $U$ is compact. 
Let $M_1, \dots, M_m$ denote the connected components of $M$ which meet $U$. 
For each $i = 1, \cdots, m$, we set $(U_i, N_i, X_i) = (U, N, X) \cap M_i$ and $\mu_i = \mu|_{M_i}$. 
By Case 1, the triple $(M_i, U_i, N_i)$ in $M_i$ has WEP for 
$(G_i, H_i) = \big( {\cal H}_{X_i}(M_i; \mu_i\reg)^{N_i}, {\cal H}_{X_i}(M_i;  \mu_i)\big)$. 
Since the pair $(G_i, H_i)$ can be canonically identified with 
the subpair $(G_X^N(M_i), H_X(M_i))$ of $(G_X^N, H_X)$ and 
${\cal E}^{G_i}(U_i, M_i) = {\cal E}^{G_X^N(M_i)}(U_i, M) = 
{\cal E}^{G_X^N}(U_i, M) \cap {\cal E}(U_i, M_i)$, which is an open subset of 
${\cal E}^{G_X^N}(U_i, M)$, it is seen that 
the triple $(M_i, U_i, N_i)$ in $M$ has WEP$(G_X^N, H_X)$. 
Hence, by Lemma~\ref{l_WEP_disjoint} $(\bigcup_i M_i, U, N)$ has WEP$(G_X^N, H_X)$ and by Lemma~\ref{l_WEP} so is $(M, U, N)$. 
\end{proof} 

\begin{theorem}\label{thm_main_ext}
Suppose $C$ is a compact subset of $M$, $U \in {\cal B}(M)$ is a neighborhood of $C$ in $M$ and $X$ is a closed subset of $\partial M$ with $X \cap C = \emptyset$. 
Then the triple $(M, U, C)$ has WEP$(G_X^{U, C}, H_X)$.  
\end{theorem} 

\begin{proof} 
By Lemma~\ref{lem_connect}\,(2)(i) and Remark~\ref{rem_bdy}, there exists a compact $n$-submanifold $N$ of $M$ such that 
\[ \mbox{$C \subset {\rm Int}_M N$, \ \ $N \subset {\rm Int}_M U - X$, \ \ 
$O - N$ is connected for each $O \in {\cal C}(M - C)$ \ \ and \ \ 
$\mu(\partial N) = 0$.} \] 
We show that $G^{U, C} \subset G^N$.  
Take any $h \in G^{U, C}$. 
Since $h|_U \in {\cal E}(U, M;\mu)$, we have $h|_N \in {\cal E}(N, M;\mu)$. 
By the choice of $N$, for each $L \in {\cal C}(M - N)$ there exists a unique $O \in {\cal C}(M - C)$ such that $L = O - N$. 
Since $h \in G^{U, C}$, we have $\mu(h(O)) = \mu(O)$.
Since $h|_U \in {\cal E}(U, M;\mu)$, $O \cap N \subset N \subset U$ 
and $N$ is compact, it follows that 
$\mu(h(O \cap N)) = \mu(O \cap N) \leq \mu(N) < \infty$. 
Hence, $\mu(h(L)) = \mu(L)$. 
This means that $h \in G^N$. 

By Lemma~\ref{lem_ext_mp} 
the triple $(M, U, N)$ has WEP$(G_X^N, H_X)$ and  
by Lemma~\ref{l_WEP} we conclude that 
the triple $(M, U, C)$ has WEP$(G_X^{U, C}, H_X)$. 
\end{proof} 

Since $H_X \subset G^{U, C}_X$, the next statement is an immediate consequence of Theorem~\ref{thm_main_ext}
and Lemma~\ref{l_WEP}.  

\begin{corollary}\label{cor_ext_mp} 
Suppose $C$ is a compact subset of $M$, $U \in {\cal B}(M)$ is a neighborhood of $C$ in $M$ and $X$ is a closed subset of $\partial M$ with $X \cap C = \emptyset$. 
Then the triple $(M, U, C)$ has WEP$({\cal H}_X(M; \mu))$.  
\end{corollary} 


\subsection{The weak extension theorem for ${\rm ker}\,c^\mu$} \mbox{} 

Suppose $M$ is a {\em connected} $n$-manifold and $\mu \in {\cal M}_g^{\partial}(M)$. 
In this section we deduce 
the weak extension theorem for the group $F = {\rm ker}\,c^\mu$ (Theorem~\ref{thm_ext_ker}). 
(Recall that $G = {\cal H}(M; \mu\reg)$ and $H = {\cal H}(M; \mu)$.
Note that $H_c  =F_c$ and $H(C) = F(C)$ for any compact subset $C$ of $M$.) 

\begin{theorem}\label{thm_ext_ker}
Suppose $C$ is a compact subset of $M$, 
$U$ and $V$ are open neighborhoods of $C$ in $M$ such that  
$V \cap O$ is connected for each $O \in {\cal C}(M - C)$. 
Then, the triple $(V, U, C)$ has WEP$({\rm ker}\,c^\mu, {\cal H}_c(M; \mu))$.  
\end{theorem} 

\begin{proof} 
(1) By Lemma~\ref{lem_connect}\,(2)(ii) and Remark~\ref{rem_bdy}, there exists a compact $n$-submanifold $N$ of $M$ such that \\[0.2mm] 
\hspace{15mm} $C \subset {\rm Int}_M N$, \ \ $N \subset V$, \ \ 
$N \cap O$ is connected for each $O \in {\cal C}(M - C)$ \ \ and \ \ 
$\mu(\partial N) = 0$. \\[0.2mm]
Note that ${\cal C}(N - C) = \{ N \cap O \mid O \in {\cal C}(M - C) \}$.   
Take compact subsets $D$ and $W$ of $M$ such that $C \subset {\rm Int}_M D$, $D \subset {\rm Int}_M W$ and 
$W \subset U \cap {\rm Int}_M N$. 
Since $N \subset V$ and $W \subset U$, by Lemma~\ref{l_WEP} 
it suffices to show that the triple $(N, W, C)$ has 
WEP$({\rm ker}\,c^\mu, {\cal H}_c(M; \mu))$. 

Since ${\cal E}^F(W, M) \subset {\cal E}^\ast(W, M;\mu\reg)$, 
by Theorem~\ref{thm_ext_bireg} 
there exists a neighborhood ${\cal U}$ of $i_W$ in \vspace{1.5mm} 
${\cal E}^F(W, M)$ and a map 
\hspace{8mm} $s : {\cal U} \to G(N)$ \hspace{4mm} such that \hspace{2mm} 
$s(f)|_D = f|_D$ \hspace{2mm} and \hspace{2mm} $s(i_W) = id_M$. \\[1mm] 
Replacing ${\cal U}$ by a smaller one, we may assume that $f(W) \subset N$ ($f \in {\cal U}$).
\vskip 1mm 

(2) Consider the $n$-manifold $N$ and $\mu|_N \in {\cal M}_g^\partial(N)$. 
By Theorem~\ref{thm_main_ext} 
the triple $(N, D, C)$ 
has WEP for 
$$(G', H') = \big({\cal H}_{\partial_+ N}(N;\mu|_N\reg)^{D,C}, {\cal H}_{\partial_+ N}(N;\mu|_N)\big).$$
Let 
${\cal E}^{G'}(D, N) \ \supset \ {\cal U}' 
\begin{array}[b]{c}
\sigma'_t \\[-1.2mm] 
\lra 
\end{array} H'$ be the associated LWE map. 
Each $h' \in H'$ has a canonical extension $\psi(h') \in H(N)$ and 
this defines the canonical homeomorphism $\psi : H' \cong H(N)$. 
\vskip 1mm 

(3) We show that $s(f)|_N \in G'$ for any $f \in {\cal U}$. 
Since $s(f) \in G(N)$, we have $s(f)|_N \in {\cal H}_{\partial_+N}(N;\mu|_N\reg)$. 
Since $f \in {\cal E}^F(W, M)$, 
there exists $h \in F$ such that $f = h|_W$. 
Since $s(f)|_D = f|_D = h|_D \in {\cal E}(D, M; \mu)$ and $s(f)(N) = N$, 
it follows that $s(f)|_D \in {\cal E}(D, N; \mu|_N)$. 
Take any $L \in {\cal C}(N - C)$. 
Then there exists a unique $O \in {\cal C}(M - C)$ with $L = N \cap O$. 
Let $K = O - L = O - N$. 
Consider $g \equiv h^{-1}s(f) \in {\cal H}_D(M)$. 
Since $M$ is connected, we have  
$O \cap D \neq \emptyset$ and since 
$g = id$ on $D$, we have $g(O) = O$ and so $s(f)(O) = h(O)$. 
Since $s(f) \in G(N)$, it follows that \\[0.2mm] 
\hspace{15mm} $s(f)(K) = K$ \hspace{4mm} and \hspace{4mm} 
$s(f)(L) = s(f)(O - K) = s(f)(O) - K = h(O) - K$. \\[0.2mm] 
Thus, we have $\mu(s(f)(L)) = \mu(h(O) - K)$.  
Since \\[0.2mm] 
\hspace{15mm} ${\rm Fr}_M O \subset C$, \hspace{4mm} $O - K =  L \subset N$ \hspace{4mm} and \hspace{4mm} 
$K = s(f)(K) \subset s(f)(O) = h(O)$, \\[0.2mm] 
it follows that $O \in {\cal B}_c(M)$, $K \subset O \cap h(O)$ and $O - K$ is relatively compact in $M$. 
Since $h \in F$, by Lemma~\ref{lem_ker} we have $\mu(h(O) - K) = \mu(O - K) = \mu(L)$. 
Therefore, we have $\mu(s(f)(L)) = \mu(L)$. 
\vspace{1mm} 
This means that $s(f)|_N \in G'$. 
\vskip 1mm 

(4) By (3), for any $f \in {\cal U}$, we have $s(f)|_N \in G'$ and 
$f|_D = s(f)|_D = (s(f)|_N)|_D \in {\cal E}^{G'}(D, N)$. 
Thus, we obtain the continuous map 
$\phi : {\cal U} \to {\cal E}^{G'}(D, N)$ defined by $\phi(f) = f|_D$. 
Replacing ${\cal U}$ by a smaller one, we may assume that $\phi({\cal U}) \subset {\cal U}'$. 
Finally, the associated LWE map $S_t : {\cal U} \to H(N)$ for WEP$(F, H_c)$ of the triple $(N, W, C)$ 
is defined by 
$$S_t(f) = \psi \sigma'_t \phi(f).$$  
\vskip -6mm
\end{proof} 

Since $H_c \subset F$, the next statement is an immediate consequence of Theorem~\ref{thm_ext_ker} 
and Lemma~\ref{l_WEP}.  

\begin{corollary}\label{cor_ext_cpt}
Suppose $C$ is a compact subset of $M$, 
$U$ and $V$ are open neighborhoods of $C$ in $M$ such that  
$V \cap O$ is connected for each $O \in {\cal C}(M - C)$. 
Then the triple $(V, U, C)$ has WEP$({\cal H}_c(M; \mu))$.  
\end{corollary} 

\subsection{Non-ambient weak deformation of measure-preserving embeddings} \mbox{} 

Suppose $M$ is an $n$-manifold and $\mu \in {\cal M}_g^\partial(M)$. 
In this section we obtain 
a non-ambient weak deformation theorem for measure-preserving embeddings.  
For $X \in {\cal B}(M)$, let 
${\cal E}^\ast(X, M; \mu) = {\cal E}(X, M; \mu) \cap {\cal E}^\ast(X, M)$ with the compact-open topology. 

\begin{theorem}\label{thm_non-amb_def} 
Suppose $C$ is a compact subset of $M$ and $U \in {\cal B}(M)$ is a neighborhood of $C$ in $M$. 
Then there exists a neighborhood ${\cal U}$ of $i_U$ in ${\cal E}^\ast(U, M; \mu)$ 
and a map $s : {\cal U} \times [0,1] \to {\cal E}^\ast(C, M; \mu)$ such that 
$s_0(f) = i_C$, 
$s_1(f) = f|_C$ $(f \in {\cal U})$ and 
$s_t(i_U) = i_C$ $(t \in [0,1])$. 
\end{theorem}  

We call the map $s$ a {\em local weak deformation} map (a LWD map) for the pair $(U, C)$ in $M$. 

\begin{lemma}\label{lem_non-amb_def} 
Suppose $N$ is a compact $n$-submanifold of $M$ with $\mu(\partial_+ N) = 0$ and $U \in {\cal B}(M)$ is a neighborhood of $N$ in $M$. 
Then the pair $(U, N)$ admits a LWD map in $M$. 
\end{lemma}  

\begin{proof}
{Case 1:}  First we treat the case where $N$ is connected. 

(1) By Lemma~\ref{lem_collar} 
there exists a bicollar $E = \partial_+ N \times [a, b]$ ($a < 0 < b$) 
of $\partial_+ N$ in $M$  
such that 
\[ \mbox{$\partial_+ N = \partial_+ N \times \{ 0 \}$, \ \ $N \cap E = \partial_+ N \times [a, 0]$ \hspace{2mm} 
and \hspace{2mm} $\mu|_E = \nu \times (m|_{[a, b]})$,} \]  
where $\nu \in {\cal M}_g^\partial(\partial_+ N)$ and $m$ is the Lebesgue measure on ${\Bbb R}$.  
Let ${\cal C}(\partial_+ N) = \{ F_1, \cdots, F_m \}$ and 
$E_i = F_i \times [a, b]$ ($i = 1, \cdots, m$). 
For notational simplicity, we use the following notations:  \ \ 
\[ \mbox{$E(I) = \partial_+ N \times I$, \ \ $E_i(I) = F_i \times I$ \ \ ($I \subset [a,b]$) \hspace{2mm} 
and \hspace{2mm} $N_t = (N - E) \cup E[a,t]$ \ \ ($t \in [a,b]$).} \] 
Take $\e > 0$ such that $a < -3 \e$, $3 \e < b$, and 
define $\alpha_t \in {\cal H}_\partial([a,b])$ ($t \in (-2\e, 2\e)$) by the conditions: \\
\hspace{15mm} $\alpha_t(s) = s + t$ ($s \in [-\e,\e]$) \ \ and \ \ $\alpha_t$ is affine on the intervals $[a , -\e]$ and $[\e, b]$. \\[0.5mm] 
For each $i = 1, \cdots, m$, we obtain the isotopy \hspace{2mm}  
$\phi^i_t = id_{F_i} \times \alpha_t \in {\cal H}_{\partial_+ E_i}(E_i; \mu|_{E_i}\reg)$ \ \ ($t \in (-2\e, 2\e)$). \\ 
Note that $\alpha_0 = id_{[a,b]}$ and $\phi^i_0 = id_{E_i}$. 

Take a small neighborhood ${\cal W}$ of $i_N$ in  ${\cal E}^\ast(N, M; \mu\reg)$ such that 
for any $g \in {\cal W}$ and $i = 1, \cdots, m$, 
\[ \mbox{$E_i[a,-\e] \subset g(N) \cap E_i \subset E_i[a, \e]$, \ \ 
$N_{-\e} \subset g(N) \subset N_\e$ \hspace{2mm} and \hspace{2mm} 
$g(F_i) \subset E_i(-\e, \e)$.} \]  
Then, for each $g \in {\cal W}$ and $i = 1, \cdots, m$, we have  
\begin{itemize}
\item[(i)\,] $(-\e - a)\nu(F_i) < \mu(g(N) \cap E_i) < (\e - a)\nu(F_i)$, 
\item[(ii)] $\mu(\phi_t^i(g(N) \cap E_i)) = \mu(g(N) \cap E_i) + t \nu(F_i)$, \  
since $\phi_t^i$ is $\mu$-preserving on $E_i[-\e, \e]$. 
\end{itemize}
For each $i = 1, \cdots, m$, consider the map \hspace{4mm}  
$c_i : {\cal W} \to {\Bbb R}$ \ \text{defined by} \ 
$c_i(g) = \mu(g(N) \cap E_i)$.  \\ 
Since $\mu(g(\partial_+ N)) = 0$, the map $c_i$ is seen to be continuous. 
Note that $c_i(g)\in \big((-\e - a)\nu(F_i), (\e - a)\nu(F_i)\big)$. 

\vskip 2mm 

(2) Next we construct a neighborhood ${\cal U}$ of $i_U$ in ${\cal E}^\ast(U, M; \mu)$ and 
a map $\eta : {\cal U} \times [0,1] \to {\cal E}^\ast(N, M; \mu\reg)$ such that 
for any $f \in {\cal U}$ and $t \in [0,1]$,  
\[ \mbox{(iii) \ $\eta_0(f) = i_N$, \ $\eta_1(f) = f|_N$, \  
$\eta_t(i_U) = i_N$ \hspace{2mm} and \hspace{2mm} (iv) \ 
$\mu(\eta_t(f)(N)) = \mu(N)$.} \] 

By Theorem~\ref{thm_ext_bireg} there exists a neighborhood ${\cal U}$ of $i_U$ in ${\cal E}^\ast(U, M; \mu)$ and 
a map 
\[ \mbox{$\sigma : {\cal U} \times [0,1] \to {\cal H}_c(M;\mu\reg)$ such that 
$\sigma_0(f) = id_M$,
$\sigma_1(f)|_N = f|_N$ ($f \in {\cal U}$) and 
$\sigma_t(i_U) = id_M$ ($t \in [0,1]$).} \] 
Replacing ${\cal U}$ by a smaller one, we may assume that   
$\sigma_t(f)|_N \in {\cal W}$ ($f \in {\cal U}$, $t \in [0,1]$). 
Consider the map 
\[ \mbox{$\gamma : {\cal U} \times [0,1] \to {\cal W} \ \subset \ {\cal E}^\ast(N, M;\mu\reg)$ \hspace{4mm} \text{defined by} \hspace{4mm}$\gamma_t(f) = \sigma_t(f)|_N$.} \]  
The map $\gamma$ satisfies the condition (iii). To achieve the condition (iv) we modify the map $\gamma$. 

We define the maps \ $\lambda^i : {\cal U} \times [0,1] \to {\Bbb R}$ \ and \ $\tau^i : {\cal U} \times [0,1] \to (-2\e, 2\e)$ \ by 
\[ \mbox{$\lambda_t^i(f) = (1-t)c_i(i_N) + tc_i(f|_N)$ 
\hspace{2mm} and \hspace{2mm} 
$c_i(\gamma_t(f)) + \tau_t^i(f) \nu(F_i) = \lambda_t^i(f)$.} \] 
Since $\lambda_t^i(f),  c_i(\gamma_t(f)) \in \big((-\e - a)\nu(F_i), (\e - a)\nu(F_i)\big)$, 
we have \   
$$|\tau_t^i(f)| \nu(F_i) = |\lambda_t^i(f) - c_i(\gamma_t(f))| 
< 2\e\nu(F_i).$$ 
The map $\tau^i$ has the following properties:  
\begin{itemize}
\item[(v)] $\tau^i_0(f) = \tau^i_1(f) = \tau^i_t(i_U) = 0$, \hspace{2mm} 
\item[(vi)] $\mu(\phi^i_{\tau_t^i(f)}(\gamma_t(f)(N) \cap E_i)) 
= \mu(\gamma_t(f)(N) \cap E_i) + \tau_t^i(f) \nu(F_i)
= \lambda_t^i(f)$.
\end{itemize} 
The assertion (vi) follows from the property (1)(ii), while 
the assertion (v) follows from \\ 
\hspace{1mm} 
\begin{tabular}[t]{l} 
$\tau^i_0(f)\nu(F_i) = \lambda^i_0(f) - c_i(\gamma_0(f)) 
= c_i(i_N) - c_i(i_N) = 0$, \hspace{1mm} 
$\tau^i_1(f)\nu(F_i) = \lambda^i_1(f) - c_i(\gamma_1(f)) 
= c_i(f|_N) - c_i(f|_N) = 0$,  \\[2mm] 
$\tau^i_t(i_U)\nu(F_i) = \lambda^i_t(i_U) - c_i(\gamma_t(i_U)) 
= c_i(i_N) - c_i(i_N) = 0.$ 
\end{tabular} \\[1.5mm] 
The maps $\phi^i_{\tau^i}$ $(i = 1, \cdots, m)$ are combined to induce the map  
\[ \mbox{$\phi : {\cal U} \times [0,1] \to {\cal H}_{E^c}(M; \mu\reg)$ \ \text{defined by} \ 
$\phi_t(f)|_{E_i} = \phi^i_{\tau_t^i(f)}$ \hspace{2mm} $(i = 1, \cdots, m)$.} \] 
The desired map \hspace{5mm} 
$\eta : {\cal U} \times [0,1] \to {\cal E}^\ast(N, M; \mu\reg)$ \hspace{4mm} \text{is defined by} \hspace{4mm} $\eta_t(f) = \phi_t(f)\gamma_t(f)$. \\
From (v) it follows that 
$\phi_0(f) = \phi_1(f) = \phi_t(i_U) = id_M$, \ 
since \ 
$$\phi_0(f)|_{E_i} = \phi_1(f)|_{E_i} = \phi_t(i_U)|_{E_i} = \phi_0^i = id_{E_i}.$$ 
Thus, the map $\eta$ satisfies the condition (iii). 
To see the condition (iv), first note that 
$$\eta_t(f)(N) = \phi_t(f)\gamma_t(f)(N) 
= \phi_t(f) \Big(N_a \cup \big(\mbox{$\bigcup_i$} (\gamma_t(f)(N) \cap E_i) \big)\Big)
= N_a \cup \Big(\mbox{$\bigcup_i$}\,\phi^i_{\tau_t^i(f)}(\gamma_t(f)(N) \cap E_i) \Big).$$
Since $f$ is $\mu$-preserving, we have $\mu(f(N)) = \mu(N)$. 
Therefore, from (vi) it follows that 
$$\begin{array}[t]{lll}
\mu(\eta_t(f)(N)) 
&=& \mu(N_a) + \sum_i \mu(\phi^i_{\tau_t^i(f)}(\gamma_t(f)(N) \cap E_i)) 
\ = \ \mu(N_a) + \sum_i \lambda_t^i(f) \\[3mm] 
&=& \mu(N_a) + (1-t) \sum_i c_i(i_N) + t \sum_i c_i(f|_N) \\[2mm] 
&=& (1-t) \Big( \mu(N_a) + \sum_i c_i(i_N) \Big) + t \Big( \mu(N_a) + \sum_i c_i(f|_N) \Big) \\[3mm] 
&=& (1-t) \mu(N) + t \mu(f(N)) \ = \ \mu(N). 
\end{array}$$
\vskip 2mm 
(3) The required LWD map $s$ is obtained as follows. 

Theorem~\ref{thm_selection} yields a map \ \ $\chi : {\cal M}_g^\partial(N; \mu|_N) \to {\cal H}_\partial(N; \mu|_N\reg)_1$ \ \ 
such that 
\[ \mbox{ $\chi (\omega)_\ast (\mu|_N) = \omega$ \ $(\omega \in {\cal M}_g^\partial(N; \mu|_N))$ \hspace{2mm} and \hspace{2mm} $\chi (\mu|_N) = id_N$.} \] 
By the condition (2)(iv) we have the map \hspace{2mm} 
$\rho : {\cal U} \times [0,1] \to {\cal M}_g^\partial(N;\mu|_N)$ \ \text{defined by} \ 
$\rho_t(f) = \eta_t(f)^\ast \mu$. \\
Since $\rho_t(f) = \eta_t(f)^\ast \mu = \big((\phi_t(f) \sigma_t(f))^\ast \mu \big)|_N$, the map $\rho$ is the composition of the following maps: 
\[ \begin{array}[c]{l} 
\begin{array}[c]{ccccccc} 
& \rho_1 & & \rho_2 & & \rho_3 & \\ 
{\cal U} \times [0,1] & \lra & {\cal H}(M;\mu\reg) & \lra & {\cal M}_g^\partial(M;\mu\reg) & 
\lra & {\cal M}_g^\partial(N;\mu|_N\reg), \\[2mm]  
\end{array} \\[4mm] 
\begin{array}[c]{c} 
\hspace{5mm} \text{where} \hspace{6mm} \rho_1(f, t) =\phi_t(f) \sigma_t(f), \hspace{6mm} 
\rho_2(h) = h^\ast \mu \hspace{6mm} \text{and} \hspace{6mm} \rho_3(\omega) = \omega|_N. 
\end{array}
\end{array} \] 
Since $\mu(\partial_+ N) = 0$, by Lemma~\ref{lem_rest} the third map is continuous. 
Thus the continuity of the map $\rho$ follows from the continuity of these maps.  
Finally, the map 
\[ \mbox{$s : {\cal U} \times [0,1] \to {\cal E}^\ast(N, M; \mu)$ \hspace{4mm} \text{ is defined by } \hspace{4mm}$s_t(f) = \eta_t(f) \chi(\rho_t(f))$.} \] 
Since $s_t(f)^\ast \mu = \chi(\rho_t(f))^\ast (\eta_t(f)^\ast \mu)
= \chi(\rho_t(f))^\ast\rho_t(f) = \mu|_N$, it follows that 
$s_t(f)$ is $\mu$-preserving. 
If $t = 0, 1$ or $f = i_U$, then by (2)(iii), $\eta_t(f)$ is $\mu$-preserving, and so $\rho_t(f) = \mu|_N$ and $s_t(f) = \eta_t(f)$. 
Hence, by (2)(iii) the map $s$ satisfies the required conditions: 
$s_0(f) = i_N$, $s_1(f) = f|_N$ and $s_t(i_U) = i_N$. 
\smallskip 

\noindent {Case 2:} Next we treat the general case where $N$ may not be connected. 

Let ${\cal C}(N) = \{ N_1, \cdots, N_m \}$. 
By Case 1, 
each pair $(U, N_i)$ ($i = 1, \cdots, m$) admits a LWD map in $M$ 
$${\cal E}^\ast(U, M; \mu) \ \supset \ {\cal U}_i 
\begin{array}[b]{c}
s^i_t \\[-0.5mm] 
\lra 
\end{array}
{\cal E}^\ast(N_i, M; \mu) \ \ (t \in [0,1]).$$ 
For each $i = 1, \cdots, m$, choose a neighborhood $U_i$ of $N_i$ in $U$ such that $U_i \cap U_j = \emptyset$ ($i \neq j$). 

We can find a small neighborhood 
${\cal U}$ of $i_U$ in ${\cal E}^\ast(U, M; \mu)$ 
such  that ${\cal U} \subset {\cal U}_i$ and 
$s^i_t(f)(N_i) \subset U_i$ ($f \in {\cal U}$) for each $i = 1, \cdots, m$. 
A LWD map 
\[ \mbox{$s : {\cal U} \times [0,1] \to {\cal E}^\ast(N, M; \mu)$ \hspace{4mm} for $(U, N)$ is defined by \hspace{4mm} 
$s_t(f)|_{N_i} = s^i_t(f)$ \ \ ($i = 1, \cdots, m$).} \] 
\vskip -7mm 
\end{proof} 

\begin{proof}[\bf Proof of Theorem~\ref{thm_non-amb_def}]
By Lemma~\ref{lem_submfd} and Remark~\ref{rem_bdy} there exists 
a compact $n$-submanifold $N$ of $M$ 
such that $\mu(\partial_+ N) = 0$ and $C \subset N \subset {\rm Int}_M U$. 
By Lemma~\ref{lem_non-amb_def} the pair $(U, N)$ admits a LWD map 
$${\cal E}^\ast(U, M; \mu) \ \supset \ {\cal U} 
\begin{array}[b]{c}
\sigma_t \\[-0.5mm] 
\lra 
\end{array}
{\cal E}^\ast(N, M; \mu) \ \ (t \in [0,1]).$$ 
A LWD map \hspace{4mm} $s_t : {\cal U} \to {\cal E}^\ast(C, M; \mu)$ \hspace{4mm} for $(U, C)$ is defined by \hspace{4mm} 
$s_t(f) = \sigma_t(f)|_C$. 
\end{proof} 


\section{Groups of measure preserving homeomorphisms endowed with the Whitney topology} 

Suppose $M$ is a {\em connected noncompact} $n$-manifold and $\mu \in {\cal M}_g^{\partial}(M)$. 
In \cite[Proposition 5.3]{BMSY} we have shown that 
the group ${\cal H}_c(M)_w$, endowed with the Whitney topology,  is locally contractible. 
In this section we shall apply the weak extension theorem for ${\cal H}_c(M; \mu)$ (Corollary~\ref{cor_ext_cpt}) to verify the local contractibility of the group 
${\cal H}_c(M; \mu)_w$ endowed with the Whitney topology (Theorem~\ref{thm_LC}). 

\subsection{Homeomorphism groups with the Whitney topology} \mbox{} 

First we recall basic properties of the Whitney topology on homeomorphism groups (cf. \cite[Section 4.3]{BMSY}). 
Suppose $Y$ is a paracompact space and 
$\cov(Y)$ is the family of all open covers of $Y$. 
For maps $f,g : X \to Y$ and $\U\in\cov(Y)$, we say that 
$f,g$ are {\em $\U$-near} and write $(f,g)\prec\U$ 
 if every point $x\in X$ admits $U\in\U$ 
 with $f(x),g(x) \in U$. 
For each $h \in {\cal H}(Y)$ and $\U\in\cov(Y)$, let \hspace{2mm} 
$$\U(h)=\{f\in {\cal H}(Y) \mid (f,h)\prec\U\}.$$
The Whitney topology on ${\cal H}(Y)$ is generated by 
the base $\U(h)$ ($h \in {\cal H}(Y)$, $\U \in \cov(Y)$). 
The symbol ${\cal H}(Y)_w$ denotes 
the group ${\cal H}(Y)$ endowed with the Whitney topology 
(while the symbol ${\cal H}(Y)$ denotes the group ${\cal H}(Y)$ with the compact-open topology). 
It is known that $G = {\cal H}(Y)_w$ is a topological group. 
Recall the notations $G_0 = {\cal H}_0(Y)_w$ (the identity component of $G$) and $G_c = {\cal H}_c(Y)_w$ (the subgroup of $G$ consisting of homeomorphisms with compact support). 
In \cite[Sections 4.1, 4.3]{BMSY} it is shown that ${\cal H}_0(Y)_w \subset {\cal H}_c(Y)_w$. 

\subsection{The box topology on topological groups} \mbox{} 

The Whitney topology is closely related to box products (cf. \cite{BMSY}). 
Next we recall basic properties of (small) box products (cf. \cite[Sections 1, 2]{BMSY}). 
The {\em box product} \ $\square_{n \geq 1}X_n$
 of a sequence of topological spaces $(X_n)_{n \geq 1}$ 
is the product $\prod_{n \geq 1}X_n$
 endowed with the box topology generated 
 by the base consisting of boxes $\prod_{n \geq 1}U_n$ 
 ($U_n$ is an open subset of $X_n$). 
The {\em small box product} \ $\cbox_{n \geq 1}X_n$ of 
a sequence of pointed spaces $\big((X_n,*_n)\big)_{n \geq 1}$  is 
the subspace of $\square_{n \geq 1}X_n$ defined by 
$$\cbox_{n \geq 1}X_n = \big\{(x_n)_{n \geq 1}\in\square_{n \geq 1}X_n \mid 
 \exists \,m\geq 1 \ \text{such that} \  x_n=*_n \ (n\ge m) \big\}.$$ 
It has the canonical distinguished point $(\ast_n)_{n \geq 1}$. 
For a sequence of subsets $A_n \subset X_n$ ($n \geq 1$), 
we set 
$$\cbox_{n \geq 1}A_n = \cbox_{n \geq 1}X_n \cap \square_{n \geq 1}A_n.$$ 

We say that a space $X$ is {\em $($strongly$)$ locally contractible} at $x \in X$
 if every neighborhood $V$ of $x$ contains
 a neighborhood $U$ of $x$ which is contractible in $V$ (rel. $x$) 
 (i.e., there is a homotopy $h : U \times [0,1] \to V$ 
 such that $h_0 = id_U$, $h_1(U) = \{ x \}$ (and 
$h_t(x) = x$ ($t \in [0,1]$)). 
A pointed space $(X, x_0)$ is said to be {\em locally contractible} 
if $X$ is locally contractible at any point of $X$ and 
strongly locally contractible at $x_0$.
It is easily seen that if a topological group $G$ is locally contractible at the identity element $e$,  
then the pointed space $(G,e)$ is locally contractible  (\cite[Remark 1.9]{BMSY}). 
The next lemma follows from a straightforward argument. 

\begin{lemma}\label{lem_LC}{\rm (\cite[Proposition 1.10]{BMSY})} 
 If pointed spaces $(X_i, \ast_i)$ $(i \geq 1)$ are locally contractible,
 then the small box product $\cbox_{i \geq 1} (X_i, \ast_i)$ 
 is also locally contractible as a pointed space. 
\end{lemma}

Suppose $G$ is a topological group with the identity element $e \in G$.    
A sequence of closed subgroups $(G_n)_{n \geq 1}$ of $G$ is called 
a {\em tower} in $G$ if it satisfies the following conditions: 
\[ G_1 \subset G_2 \subset G_3 \subset \cdots \hspace{5mm} \text{and} \hspace{5mm} G = \mbox{$\bigcup_{n \geq 1}$} G_n. \] 
Any tower $(G_n)_{n \geq 1}$ in $G$ induces 
the small box product $\cbox_{n \geq 1}(G_n, e)$ and 
the multiplication map 
\[ p :\cbox_{n \geq 1}(G_n, e) \lra G \hspace{5mm} \text{defined by} \hspace{5mm} 
p(x_1, \dots, x_m, e, e, \cdots) = x_{1}\cdots x_m. \]  
Note that $\cbox_{n \geq 1} G_n$ is 
a topological group with the coordinatewise multiplication and 
the  identity element $\boldsymbol{e} = (e, e, \cdots)$ and that 
the map $p$ is well-defined and continuous (\cite[Lemma 2.1]{BMSY}). 

\begin{defi}\label{def_box_top}
We say that 
$G$ {\em carries the box topology}
 with respect to $(G_n)_{n \geq 1}$ if the map $p:\cbox_{n \geq 1}G_n\to G$ is an open map. 
\end{defi}

Recall that $G$ is {\em the direct limit} of $(G_n)_{n \geq 1}$ 
in the category of topological groups 
 if any group homomorphism $h : G \to H$ to an arbitrary topological group $H$
 is continuous provided the restriction $h|G_n$ is continuous 
 for each $n \geq 1$. 
If $G$ carries the box topology with respect to $(G_n)_{n \geq 1}$, then 
$G$ is the direct limit of $(G_n)_{n \geq 1}$ in the category of topological groups (\cite[Proposition 2.7]{BMSY}). 
Note that the map $p$ is an open map if it is open at $\boldsymbol{e}$ 
(i.e., for any neighborhood $U$ of $\boldsymbol{e}$ in $\cbox_{n \geq 1}G_n$ 
 the image $p(U)$ is a neighborhood of $e$ in $G$).  
We say that a map $f : X \to Y$ has a local section at $y \in Y$ if there exists a neighborhood $U$ of $y$ in $Y$ and a map 
$s : U \to X$ such that $fs = i_U$. 
If the map $p$ has a local section $s : U \to \cbox_{n \geq 1} G_n$ at $e \in G$, then 
(i) we can adjust $s$ so that $s(e) = \bs{e}$ and so 
(ii) the map $p$ is open at $\bs{e}$. 
Thus, the next lemma follows from Definition~\ref{def_box_top} and Lemma~\ref{lem_LC}. 

\begin{lemma}\label{lem_sec_LC} Suppose the map $p:\cbox_{n \geq 1}G_n\to G$ has a local section at $e$. Then 
\begin{itemize}
\item[(1)] $G$ carries the box topology 
 with respect to the tower $(G_n)_{n \geq 1}$, 
\item[(2)] if the subgroups $G_n$ $(n \geq 1)$ are locally contractible, 
then $G$ is also locally contractible. 
\end{itemize}
\end{lemma} 

\begin{lemma}\label{lem_loc_sec_subseq}
The map $p:\cbox_{n \geq 1}G_n\to G$ has a local section at $e$ iff  
for any $($or some\,$)$ subsequence $(G_{n(i)})_{i \geq 1}$
the multiplication map 
$p' :\cbox_{i \geq 1}G_{n(i)} \lra G$ has a local section at $e$. 
 \end{lemma}

\begin{proof} 
Consider the maps 
\hspace{4mm} $\mbox{$\pi : \cbox_{n \geq 1}G_n \to \cbox_{i \geq 1}G_{n(i)}$ \hspace{3mm} and  \hspace{3mm} 
$\eta : \cbox_{i \geq 1}G_{n(i)} \to \cbox_{n \geq 1}G_n$}$ 
\vspace{-2mm} 
\[ \mbox{defined by \hspace{4mm} 
$\begin{array}[t]{l}
\pi( \cdots, x_{n(i-1)+1}, \cdots, x_{n(i)}, \cdots) = 
( \cdots, \!\!\!
\begin{array}[b]{c}
\mbox{\tiny $i$} \\[-2mm] 
\mbox{\small $\vee$} \\[-1mm] 
(x_{n(i-1)+1} \cdots x_{n(i)})
\end{array} \!\!, \cdots) \hspace{5mm} \text{and} \\[3mm]  
\eta(\cdots, x_{i-1}, x_{i}, \cdots ) = 
(\cdots, e, 
\hspace{-1.5mm} 
\begin{array}[t]{c}
x_{i-1} \\[-1mm] 
\mbox{\small $\wedge$} \\[-2mm] 
\mbox{\tiny $n(i\!-\!1)$} 
\end{array}\hspace{-1.5mm} , e, \cdots, e, \hspace{-1.5mm} 
\begin{array}[t]{c}
x_i \\[-1mm] 
\mbox{\small $\wedge$} \\[-2mm] 
\mbox{\tiny $n(i)$} 
\end{array}\hspace{-1.5mm} , \cdots), \hspace{5mm} \text{where $n(0) = 0$}. 
\end{array}$} \] 
\noindent 
The maps $p$ and $p'$ have the factorizations 
$p' = p\eta$ and 
$p = p'\pi$, from which follows the assertion. 
\end{proof} 

\subsection{Local contractibility of ${\cal H}_c(M; \mu)_w$} \mbox{} 

Suppose $M$ is a connected noncompact $n$-manifold and $\mu \in {\cal M}_g^\partial(M)$. Let $H = {\cal H}(M; \mu)$ and $F = {\rm ker}\,c^\mu$. 
(Recall that the subscript $w$ means the Whitney topology. 
For example, $H_{c,w} = {\cal H}_c(M; \mu)_w$.) 

Consider any sequence $(K_i)_{i \geq 1}$ of compact subsets of $M$ such that 
$K_ i \subset {\rm Int}_M K_{i+1}$ $(i \geq 1)$  
and $M=\bigcup_{i \geq 1}K_i.$ 
It induces a tower $H(K_i) = {\cal H}_{M - K_i}(M; \mu)$ $(i \geq 1)$ of $H_{c,w}$ and the multiplication map
\[ p : \cbox_{i \geq 1} H(K_i) \lra H_{c,w}, 
\quad p(h_1, \dots, h_m, id_M, id_M, \cdots) = h_1 \cdots h_m. \]

\begin{theorem}\label{thm_LC} 
{\rm (1)} The multiplication map 
$p : \cbox_{i \geq 1} H(K_i) \to {\cal H}_c(M; \mu)_w$ has a local section at $id_M$. 
\begin{itemize}
\item[(2)] The group ${\cal H}_c(M; \mu)_w$ carries the box topology 
 with respect to the tower $(H(K_i))_{i \geq 1}$. 
\item[(3)] The group ${\cal H}_c(M; \mu)_w$ is locally contractible. 
\end{itemize} 
\end{theorem}

We need some preliminary lemmas. 
Consider a sequence of compact connected $n$-submanifolds $(M_i)_{i \geq 1}$ of $M$ such that 
$M_ i \subset {\rm Int}_M M_{i+1}$ ($i \geq 1$) 
and $M = \bigcup_{i \geq 1}M_i$. 
Let $M_0 = \emptyset$ and $L_i = M_i - {\rm Int}_M M_{i-1}$ ($i \geq 1$). 
There exists a sequence of compact $n$-submanifolds $(N_i)_{i \geq 1}$ of $M$ such that 
$L_i \subset {\rm Int}_M N_i$ and $N_i \cap N_j \neq \emptyset$ 
 iff $|i - j| \leq 1$. 
We call the sequence $(M_i, L_i, N_i)_{i \geq 1}$ an {\em exhausting sequence} for $M$. 

\begin{lemma}\label{lem_exh_seq} 
For any sequence $(K_i)_{i \geq 1}$ of compact subsets of $M$ 
there exists an exhausting sequence $(M_i, L_i, N_i)_{i \geq 1}$ for $M$ such that 
for each $i \geq 1$ {\rm (i)} $K_i \subset M_i$, {\rm (ii)} $\mu(\partial_+ M_i) = 0$ and {\rm (iii)} the pair $(N_i, L_i)$ has WEP$(F, H_c)$. 
\end{lemma} 

\begin{proof}
By the repeated application of Lemma~\ref{lem_submfd},  
we can find a sequence of compact connected $n$-submanifolds 
$(M_i)_{i \geq 1}$ of $M$ such that 
\begin{itemize} 
\item[(i)\,] $K_i \subset M_i \subset {\rm Int}_M M_{i+1}$, $\mu(\partial_+ M_i) = 0$ $(i \geq 1)$ and $M = \bigcup_{i \geq 1}M_i$, 
\item[(ii)] $L$ is noncompact and $M_{i+1} \cap L$ is connected for each $i \geq 1$ and each $L \in {\cal C}(M_i^c)$. 
\end{itemize} 
Let $M_i = \emptyset$ $(i \leq 0)$ and $M_i^j = M_j - {\rm Int}_M M_i$ ($j > i$). 
\vskip 1mm 
(1) First we show that the pair $(N, K) = (M_{i-1}^{j+1}, M_i^j)$ has WEP$(F, H_c)$
for each $j > i \geq 0$. 
Let ${\cal C}(M_{i-1}^c) = \{ C_1, \cdots, C_m \}$ and set 
$(N_k, K_k) = (N \cap C_k, K \cap C_k)$ ($k = 1, \cdots, m$).  
Since $(N_k)_k$ is a disjoint finite family, by Lemma~\ref{l_WEP_disjoint} 
it suffices to show that each pair $(N_k, K_k)$ has WEP$(F, H_c)$. 

Note that ${\cal C}(K_k^c) = \{ E_0, E_1, \cdots, E_\ell\}$, where 
$$E_0 = M_i \cup \mbox{$\bigcup_{s \neq k}$} C_s \hspace{4mm} \text{and} \hspace{4mm} 
\{ E \in {\cal C}(M_j^c) \mid E \subset C_k \} = \{ E_1, \cdots, E_\ell \}. $$ 
(If $i=0$, we ignore $E_0$.) 
By the above condition (ii) it is seen that the intersections 
\[ \mbox{$N_k \cap E_0 = M_i \cap C_k$ \hspace{3mm} and \hspace{3mm}  
$N_k \cap E_t = M_{j+1} \cap E_t$ \hspace{2mm} ($t = 1, \cdots, \ell$)} \]
 are connected. 
Hence, we can apply Theorem~\ref{thm_ext_ker} to 
$(V, U, C) = ({\rm Int}_M N_k, {\rm Int}_M N_k, K_k)$ to conclude that 
this triple has WEP$(F, H_c)$. 
Thus, by Lemma~\ref{l_WEP} the pair $(N_k, K_k)$ also has WEP$(F, H_c)$. 
\vskip 1mm 

(2) Now consider the subsequence $(M_{3i})_{i \geq 1}$.  
Let $L_i = M_{3i-3}^{3i}$ and $N_ i = M_{3i-4}^{3i+1}$ ($i \geq 1$). 
Then, it is seen that $(M_{3i}, L_i, N_i)_{i \geq 1}$ is an exhausting sequence for $M$ and by (1) each pair $(N_i, L_i)$ has WEP$(F, H_c)$. 
\end{proof}

Suppose $(M_i, L_i, N_i)_{i \geq 1}$ is an exhausting sequence for $M$. 
It induces a tower $(H(M_i))_{i \geq 1}$ of $H_{c,w}$ and the multiplication map
$p : \cbox_{i \geq 1} H(M_i) \lra H_{c,w}$. 

\begin{lemma}\label{lem_loc_sec} 
If each pair $(N_{2i}, L_{2i})$ $(i \geq 1)$ has WEP$(H_c)$, then 
the map $p : \cbox_{i \geq 1} H(M_i) \to H_{c,w}$ has a local section $s : {\mathcal U} \to \cbox_{i \geq 1} H(M_i)$ at $id_M$ such that  
$s(id_M) = (id_M)_{i \geq 1}$ 
\end{lemma} 

\begin{proof} 
We use the following notations: 
Let \ $L_e = \bigcup_i  L_{2i}$, \ $L_o = \bigcup_i  L_{2i-1}$ \ and \  
$N_e = \bigcup_i  N_{2i}$. 
Consider the continuous maps defined by \\[1mm] 
\begin{tabular}[c]{cl} 
(a) & $r_e : H_{c, w} \to \cbox_i {\cal E}^{H_c}(L_{2i},M)$, \ $r_e(h) = (h|_{L_{2i}})_i $ \hspace{2mm} and \hspace{2mm} 
 $r : H_{c, w} \to \cbox_i {\cal E}^{H_c}(N_{2i},M)$, \ $r(h) = (h|_{N_{2i}})_i $, \\[1.5mm] 
(b) & $\lambda : \cbox_i  H(N_{2i}) \to {H_c}(N_e)_w$, \ 
$\lambda((g_i)_i)|_{N_{2i}} = g_i|_{N_{2i}}$ \hspace{0.5mm} and \hspace{0.5mm} \\[1.5mm] 
& $\lambda_o : \cbox_i  H(L_{2i-1}) \to H_c(L_o)_w$, \ 
$\lambda_o((h_i)_i)|_{L_{2i-1}} = h_i|_{L_{2i-1}}$, \\[1.5mm] 
(c) & $\rho : \cbox_i  H(N_{2i}) \times \cbox_i  H(L_{2i-1}) \to H_{c,w}, \ \ 
 \rho(\bs{g},\bs{h}) = \lambda(\bs{g}) \lambda_o(\bs{h}).$
\end{tabular} \\[1.5mm] 
Note that the map $\lambda_o$ is a homeomorphism, since for any $h \in H_c(L_o)$ we have 
$h = id$ on $\partial_+ M_i$ and $h(M_i) = M_i$, so that $h(L_i) = L_i$ ($i \geq 1$). 

First we construct a local section of the map $\rho$ at $id_M$. 
By the assumption, for each $i \geq 1$ 
there exists a neighborhood $\V_i$ of 
the inclusion map $i_{N_{2i}}$ in ${\cal E}^{H_c}(N_{2i}, M)$ 
and a map 
\[ \mbox{$\sigma_i : \V_i \to H(N_{2i})$ \hspace{2mm} such that \hspace{2mm} 
$\sigma_i(f)|_{L_{2i}} = f|_{L_{2i}}$ $(f \in \V_i)$ \ \ and \ \ $\sigma_i(i_{N_{2i}}) = id_M$.}\] 
Since $\cbox_i  \V_i$ is a neighborhood of $(i_{N_{2i}})_i$ in $\cbox_i {\cal E}^{H_c}(N_{2i},M)$, 
the preimage $\U=r^{-1}(\cbox_i  \V_i)$ is a neighborhood of $id_M$ in $H_{c,w}$. 
The maps $(\sigma_i)_i $ determine the continuous maps
\[ \mbox{$\sigma : \cbox_i  \V_i \lra \cbox_i  H(N_{2i}) \ \ \ \text{defined by} \ 
\ \ \sigma((f_i)_i ) = \big(\sigma_i(f_i)\big)_i $ \hspace{2mm} 
 and \hspace{2mm} 
$\eta = \lambda \,\sigma \,r : \U \longrightarrow H_c(N_e)_w$.} \]  
For each $g \in {\cal U}$ we have $\eta(g) = g$ on $L_e$ and 
$\eta(g)^{-1} g \in H_{c, L_e} = H_c(L_o)$. 
Thus we obtain the map 
 $$\mbox{$\phi : \U \to  H_c(L_o)_w$ \hspace{4mm} \text{defined by} \hspace{4mm} $\phi(g) = \eta(g)^{-1} g$.}$$  
The required local section \hspace{3mm} 
$\zeta : {\cal U} \lra \cbox_i H(N_{2i}) \times \cbox_i  H(L_{2i-1})$ \hspace{3mm} 
of the map $\rho$ is defined by 
\[ \zeta(g) = (\sigma\,r(g), \lambda_o^{-1} \phi(g)). \] 
In fact, we have 
$$\rho \zeta(g) = \rho(\sigma\,r(g), \lambda_o^{-1} \phi(g)) = \lambda(\sigma\,r(g)) \phi(g) = \eta(g) (\eta(g)^{-1} g) = g.$$  
Note that $\zeta(id_M) = ((id_M)_i, (id_M)_i)$. 

For each $h \in {\cal U}$ 
 the image $\zeta(h) = ((f_i)_i, (g_i)_i)$
 satisfies the following conditions:
\begin{itemize}
\item[(i)\ ]
$h = \lambda((f_i)_i) \, \lambda_o((g_i)_i)
 = (f_1f_2\cdots)(g_1g_2\cdots) = f_1g_1f_2g_2f_3g_3 \cdots$.
\vskip 0.5mm
\item[(ii)\,]
$f_i \in H(N_{2i}) \subset H(M_{2i+1})$, \quad
$g_i \in H(L_{2i-1}) \subset H(M_{2i-1})
 \subset H(M_{2i+2})$ \quad $(i \geq 1)$.
\vskip 0.5mm
\item[(iii)]
 $\big(id_M, id_M, f_1, g_1, f_2, g_2, \dots\big)
 \in \cbox_{i \geq 1} H(M_i)$ \quad and \quad
 $h = p\big(id_M, id_M, f_1, g_1, f_2, g_2, \dots\big)$.
\end{itemize}
Therefore, the required local section \ $s : {\cal U} \to \cbox_i H(M_i)$ \ of the map \ $p : \cbox_i H(M_i) \to H_{c, w}$ \ 
 is defined by 
\[s(h) = \big(id_M, id_M, f_1, g_1, f_2, g_2, \dots\big). \]
This completes the proof. 
\end{proof}

\begin{lemma}\label{lem_LC_cpt} Suppose $N$ is a compact $n$-manifold, $L$ is a $($locally flat\,$)$ $(n-1)$-submanifold of $\partial N$ and $\nu \in {\cal M}_g^\partial(N)$. Then 
the group ${\cal H}_L(N;\nu)$ is locally contractible. 
\end{lemma} 

\begin{proof} In \cite[Theorem 4.4]{Fa} the case where $L = \emptyset$ or $\partial N$ is verified. For the sake of completeness we include a proof. 
We may assume that $N$ is connected. 

(1) First we see that the group $G_L = {\cal H}_L(N;\nu\reg)$ is locally contractible. 
Since $G_L$ is a topological group, it suffices to show that 
it is semi-locally contractible at $id_N$, that is, a neighborhood of $id_N$ contracts in $G_L$. 
Using a collar $L \times [0,2]$ of $L$ in $N$ (cf. Lemma~\ref{lem_collar}), we have a deformation 
of $G_L$ to $G_{L \times [0,1]}$ which fixes $id_N$. Applying Theorem~\ref{thm-deform_bireg} to 
$(C, U, D, E) = (N, N, L, L \times [0,1])$, we can find 
a neighborhood of $id_N$ in 
$G_{L \times [0,1]}$ which contracts in 
$G_L$. 
These deformations are combined to yield a desired contraction of 
a neighborhood of $id_N$ in $G_L$. 

(2) Next we show that the group $H_L = {\cal H}_L(N;\nu)$ is a strong deformation retract (SDR) of $G_L$. 
By Theorem~\ref{thm_selection} the map 
$\pi : G \to {\cal M}_g^\partial(N; \nu)$ admits a section $s : 
{\cal M}_g^\partial(N; \nu) \to G_\partial \subset G_L$. 
This yields a homeomorphism of pairs 
$$H_L \times \big({\cal M}_g^\partial(N; \nu), \{ \nu \} \big) 
\approx (G_L, H_L) \ : \ (h, \omega) \longmapsto s(\omega) h.$$
Since ${\cal M}_g^\partial(N; \nu)$ admits the ``straight line contraction'' to $\{ \nu \}$, we obtain a SDR of $G_L$ onto $H_L$. 

Finally, the conclusion follows from the observations (1) and (2).  
\end{proof}

\begin{proof}[\bf Proof of Theorem~\ref{thm_LC}] 
(1), (3) By Lemma~\ref{lem_exh_seq} 
there exists an exhausting sequence $(M_i, L_i, N_i)_{i \geq 1}$ for $M$ such that 
$\mu(\partial_+ M_i) = 0$ $(i \geq 1)$ and each pair $(N_i, L_i)$ $(i \geq 1)$ has WEP$(H_c)$. 
By Lemma~\ref{lem_loc_sec}
the multiplication map $p' : \cbox_{i \geq 1} H(M_i) \to H_{c,w}$ has a local section at $id_M$. 
By Lemma~\ref{lem_loc_sec_subseq} this implies the assertion (1)  
(consider a mixed sequence of $(K_i)_i$ and $(M_i)_i$). 
By Lemma~\ref{lem_LC_cpt} the group 
 $H(M_i) \cong {\cal H}_{\partial_+ M_i}(M_i; \mu|_{M_i})$ 
 is locally contractible for each $i \geq 1$. 
 Thus, by Lemma~\ref{lem_sec_LC}\,(2) the group $H_{c,w}$ is also locally contractible. 

(2) The assertion follows from (1) and Lemma~\ref{lem_sec_LC}\,(1). 
\end{proof}


\end{document}